# ROBUST GAUSSIAN STOCHASTIC PROCESS EMULATION

By Mengyang Gu *, Xiaojing Wang † and James O. Berger ‡

*Johns Hopkins University *, University of Connecticut † and Duke University‡*

We consider estimation of the parameters of a Gaussian Stochastic Process (GaSP), in the context of emulation (approximation) of computer models for which the outcomes are real-valued scalars. The main focus is on estimation of the GaSP parameters through various generalized maximum likelihood methods, mostly involving finding posterior modes; this is because full Bayesian analysis in computer model emulation is typically prohibitively expensive.

The posterior modes that are studied arise from objective priors, such as the reference prior. These priors have been studied in the literature for the situation of an isotropic covariance function or under the assumption of separability in the design of inputs for model runs used in the GaSP construction. In this paper, we consider more general designs (e.g., a Latin Hypercube Design) with a class of commonly used anisotropic correlation functions, which can be written as a product of isotropic correlation functions, each having an unknown range parameter and a fixed roughness parameter. We discuss properties of the objective priors and marginal likelihoods for the parameters of the GaSP and establish the posterior propriety of the GaSP parameters, but our main focus is to demonstrate that certain parameterizations result in more robust estimation of the GaSP parameters than others, and that some parameterizations that are in common use should clearly be avoided. These results are applicable to many frequently used covariance functions, e.g., power exponential, Matérn, rational quadratic and spherical covariance. We also generalize the results to the GaSP model with a nugget parameter. Both theoretical and numerical evidence is presented concerning the performance of the studied procedures.

**1. Introduction.** A Gaussian Stochastic Process (GaSP) is a useful tool for analyzing spatially correlated data. For example, in geostatistics, it has been popularly used to model various types of data with complicated patterns ([10]). This paper, however, focuses on the use of GaSPs in emulation (approximation) of complex computer models. Computer models are developed in an effort to reproduce the behavior of engineering, physical,

*Keywords and phrases:* Gaussian stochastic process, emulation, anisotropic covariance, robust parameter estimation, objective priors, posterior propriety

biological and human processes. A key issue with such computer models is that they are typically very time-consuming to run (e.g., the TITAN2D computer model that models volcanic pyroclastic flows ([4]) requires up to 2 hours for a single run) and a large number of runs is typically needed for inferences concerning the computer model (i.e., estimation of parameters of the computer model) or predictions using the computer model, both being aspects of what is called *Uncertainty Quantification (UQ)* for computer models. It is thus typically crucial to develop a fast (and accurate) emulator to approximate the computer model, for use in UQ tasks ([20, 26, 3, 30]).

Data from a computer model (i.e., runs from the computer model) is typically rather different than spatial data. First, the input space of the computer model (e.g. the space of model parameters, initial conditions, boundary conditions, etc.) often has high dimension, while the maximum dimension for spatial data is typically three. Second, the inputs of a computer model typically are variables on completely different scales, so the effect of the inputs on the correlations will be highly variable. Consequently, the assumption of isotropy, which is often adopted in spatial processes, usually does not hold for modeling data from computer models. Different types of geometrically anisotropic spatial processes are discussed in the literature (c.f., [42, 18]). For computer models, it is common to use a product correlation function ([35, 28, 4, 29]), typically with very different correlation parameters for each input; the product form also keeps computations tractable, and this choice will be followed herein. Third, many computer models are deterministic, or close to being deterministic, while noise in data from spatial processes can be large. The fourth difference is that, by design, data from computer models is typically taken at input values that are far apart, whereas this may well not be so for spatial data.

In this paper, we focus on the problem of estimating the parameters of the GaSP emulator. These parameters typically consist of mean parameters, a variance parameter, and the parameters in the correlation functions, such as range and roughness parameters (introduced in more detail in the next section). Although the mean parameters and variance parameter are relatively easy to deal with, it was pointed out in [20] that the parameters in the correlation functions are notoriously difficult to estimate. For instance, maximum likelihood estimation (MLE) of these parameters has been widely recognized to be unstable ([25, 21, 23, 31]) and can be inconsistent under infill asymptotics ([40, 41]). The instability is partially caused by the Cholesky decomposition of covariance matrices that are often close to singular, when evaluating the likelihood. This can often be overcome by adding a nugget to stabilize the computation, but studies have found that the features of the emulator can significantly change when a nugget is added ([2]). Another difficulty that will be discussed herein is that serious problems can arise when the covariance matrix is estimated to be near-diagonal, and this can easily happen when a product correlation structure is used because, if even one of the terms in the product is close to zero, the correlation will be close to zero. Two R packages, DiceKriging and DiceOptim, use several different ways to avoid unstable results, such as using expected improvement criteria and bounds for the range parameters ([34]). Although these methods can yield stable computations, they produce larger predictive errors (as shown in Section 5) than the methods proposed herein, which seek parameter estimates that are naturally robust.

To obtain parameter estimates that are naturally robust, i.e., that stabilize the computation without degrading the predictive accuracy of the emulator, we utilize formal objective prior distributions (namely reference priors) and then find posterior modes for the correlation parameters. The first use of reference priors in modeling spatially correlated data was [5]; that paper was restricted to consideration of an isotropic covariance function, with only one range/scale parameter. Reference priors for an anisotropic process were studied in [28, 33], and their properties were studied in the context of product correlation functions and separable designs (e.g., a lattice) for the input values over which the computer model is run. Most designs used for creating emulators of computer models – such as the Latin Hypercube Design (LHD) – are, however, non-separable, and so we need to extend the analysis of the reference priors and likelihoods to cover non-separable situations and to include the possibility of a nugget parameter (a noise term). (Objective priors for isotropic GaSPs with a nugget were discussed recently in [6, 32, 19].)

Posterior modes of the correlation parameters depend on the parameterization used for the parameters and it was first found in [23] that this choice of parameterization can make a major difference of the "robustness" of the posterior mode. The word "robust" in this context was first used in [37] and will be formally defined in Section 3, but, informally, a robust procedure avoids the numerical issues discussed above while producing an emulator with good predictive performance. In this investigation, it was also found that robustness is considerably more difficult to obtain for the anisotropic case with product correlation functions than for the isotropic case. As an example, the posterior density of the range parameters goes to infinity when the correlation matrix, for a product correlation function, approaches a matrix of ones, under one frequently used parameterization, while this does not happen in the isotropic case. One of the major contributions of this work is in making the study of robustness of the parameterization rigorous by determining the tail behavior of the resulting posterior distributions.

The paper is organized as follows. In Section 2, we introduce the GaSP emulator with product correlation functions and designs for the input values at which the computer model is run, and we begin the comparison of our methods to two standard approaches – maximum likelihood estimation (MLE) and maximum marginal likelihood estimation (MMLE) – in order to highlight some of the key concerns. In Section 3, we first study a closed-form example of profile and marginal likelihood, where a sufficient and necessary condition is provided under which the MLE has poor behavior. Then, we formally define robust parameter estimation in the development of GaSP emulators and prove our main results concerning robustness, along with establishing posterior propriety of the suggested priors. The potentially serious consequences of using non-robust estimation methods will also be highlighted. In Section 4, we extend the results to a GaSP with a noise term. The robust method has been implemented in a new R package ([13]), which will be used for comparison of the method with other approaches, such as the MLE and DiceKriging, in Section 5. Section 6 presents some conclusions.

## 2. Gaussian stochastic processes.

2.1. *Background and a recommendation.* Consider a real-valued Gaussian stochastic process $y(\cdot) \in \mathbb{R}$ on a $p$-dimensional input domain $\mathscr{X}$,

$$y(\cdot) \sim GaSP(\mu(\cdot),\, \sigma^2 c(\cdot,\cdot)), \tag{2.1}$$

where $\mu(\cdot)$ is the mean function and $\sigma^2 c(\cdot,\cdot)$ is the covariance function with variance $\sigma^2$ and correlation function $c(\cdot,\cdot)$. For any inputs $\mathbf{x}_i \in \mathscr{X}$, $i = 1, \cdots, n$, the outputs $(y(\mathbf{x}_1), \ldots, y(\mathbf{x}_n))^T$ follow a multivariate normal distribution,

$$\left[(y(\mathbf{x}_1), \ldots, y(\mathbf{x}_n))^T \mid \boldsymbol{\mu},\, \sigma^2,\, \mathbf{R}\right] \sim \mathcal{MN}((\mu(\mathbf{x}_1), \ldots, \mu(\mathbf{x}_n))^T, \sigma^2 \mathbf{R})\,, \tag{2.2}$$

where $\mathbf{R}$ denotes the correlation matrix with the $(i,j)$ entry $c(\mathbf{x}_i, \mathbf{x}_j)$ and $\boldsymbol{\mu} = (\mu(\mathbf{x}_1), \cdots, \mu(\mathbf{x}_n))^T$. The mean function for any input $\mathbf{x} \in \mathscr{X}$ is modeled via the regression

$$\mu(\mathbf{x}) = \mathrm{E}[y(\mathbf{x})] = \mathbf{h}(\mathbf{x})\boldsymbol{\theta} = \sum_{t=1}^{q} h_t(\mathbf{x})\theta_t\,,$$

where $\mathbf{h}(\mathbf{x}) = (h_1(\mathbf{x}), h_2(\mathbf{x}), \ldots, h_q(\mathbf{x}))$ is a $q$-dimensional vector of basis functions and $\boldsymbol{\theta} = (\theta_1, \cdots, \theta_q)^T$, with $\theta_t$ being an unknown regression parameter for the basis function $h_t$.

The process is called isotropic if the correlation function is only a function of $||\mathbf{x}_i - \mathbf{x}_j||_2$, for any $\mathbf{x}_i = (x_{i1}, \cdots, x_{ip})^T \in \mathscr{X}$ and $\mathbf{x}_j = (x_{j1}, \cdots, x_{jp})^T \in \mathscr{X}$, where $|| \cdot ||_2$ is the Euclidean distance or the $L_2$ norm. As mentioned earlier, isotropy is often too restrictive to emulate complicated functions and a product of $p$ one-dimensional correlation functions is typically assumed for the computer model emulation instead

$$c(\mathbf{x}_i, \mathbf{x}_j) = \prod_{l=1}^{p} c_l(x_{il}, x_{jl}), \tag{2.3}$$

with $c_l(\cdot, \cdot)$ being a one-dimensional correlation function for the $l^{th}$ coordinate of the input vector.

The simulator is run at a set of $n$ chosen inputs $\mathbf{x}^{\mathscr{D}} = \{\mathbf{x}_1^{\mathscr{D}}, \ldots, \mathbf{x}_n^{\mathscr{D}}\}$, often selected using some "space filling" technique over the input domain $\mathscr{X}$, e.g., a Latin Hypercube Design ([35, 36]); let $\mathbf{y}^{\mathscr{D}} = (y(\mathbf{x}_1^{\mathscr{D}}), \ldots, y(\mathbf{x}_n^{\mathscr{D}}))^T$ denote the corresponding simulator outputs. Given the product correlation function in (2.3), the correlation matrix of these inputs is thus

$$\mathbf{R} = \mathbf{R}_1 \circ \mathbf{R}_2 \circ \ldots \circ \mathbf{R}_p \,, \tag{2.4}$$

where each $\mathbf{R}_l$ is the correlation matrix for the $l^{th}$ input, having $(i, j)$th element $c_l(x_{il}, x_{jl})$, and $\circ$ is the Hadamard product.

Some frequently chosen correlation functions are listed in Table 1 (dropping the subscript $l$). The correlation function $c_l(\cdot, \cdot)$ typically has a range parameter $\gamma_l > 0$, which controls how fast the correlation decays with the distance, and a roughness parameter $\alpha_l > 0$, controlling the geometric properties of the process ([5]). As mentioned earlier, the points in $\mathbf{x}^{\mathscr{D}}$ are typically chosen as far apart as possible, in order to sample the computer model output at as many diverse points as possible. Consequently, the roughness parameters $\alpha_l$, $1 \leq l \leq p$, are not highly influential and typically have quite flat likelihood surfaces. They are also highly confounded with $\gamma_l$ and $\sigma^2$, causing computational and inferential difficulties if left in the model ([40, 10]). It is thus common (and herein adopted) to fix the roughness parameters at prespecified values and focus only on estimation of the range parameters. An alternative possibility would be to assign a discrete prior – concentrated on a few values – to the roughness parameters, as in [8]; the results herein would likely generalize to that situation.

One of most frequently used correlation functions is the Gaussian correlation, which is the special case of $\alpha_l = 2$ in the power exponential correlation function. The sample paths of the resulting GaSP process are infinitely differentiable, which is sometimes desirable in applications. However, the choice

| | $c(d)$ | $\nu(\gamma)$ | $\omega(\gamma)$ |
|---|---|---|---|
| Power Exponential | $\exp\{-(d/\gamma)^\alpha\}$, $\alpha \in (0,2]$ | $\gamma^{-\alpha}$ | $\gamma^{-\alpha}$ |
| Spherical | $\left(1-\frac{3}{2}\left(\frac{d}{\gamma}\right)+\frac{1}{2}\left(\frac{d}{\gamma}\right)^3\right)\mathbf{1}_{[d/\gamma\leq 1]}$ | $\gamma^{-1}$ | $\gamma^{-2}$ |
| Rational Quadratic | $\left(1+\left(\frac{d}{\gamma}\right)^2\right)^{-\alpha}$, $\alpha \in (0,+\infty)$ | $\gamma^{-2}$ | $\gamma^{-2}$ |
| Matérn | $\frac{1}{2^{\alpha-1}\Gamma(\alpha)}\left(\frac{d}{\gamma}\right)^\alpha \mathcal{K}_\alpha\left(\frac{d}{\gamma}\right)$, $0<\alpha<1$ | $\gamma^{-2\alpha}$ | $\gamma^{-2+2\alpha}$ |
| | $\alpha = 1$ | $\frac{\log(\gamma)}{\gamma^2}$ | $\frac{1}{\log(\gamma)}$ |
| | $1<\alpha<2$ | $\gamma^{-2}$ | $\gamma^{2-2\alpha}$ |
| | $\alpha = 2$ | $\gamma^{-2}$ | $\frac{\log(\gamma)}{\gamma^2}$ |
| | $\alpha > 2$ | $\gamma^{-2}$ | $\gamma^{-2}$ |

TABLE 1

*Popular choices of correlation functions, where $c_l(x_{il}, x_{jl}) \equiv c(d)$, with $d = |x_{il} - x_{jl}|$. Here $\alpha$ is the roughness parameter and $\gamma$ is the range parameter. $\Gamma(\cdot)$ is the gamma function and $\mathcal{K}_\alpha(\cdot)$ is the modified Bessel function of the second kind. $\nu(\gamma)$ and $\omega(\gamma)$ are terms in the Taylor expansion of the correlation functions, as $\gamma \to \infty$ (see Section 3).*

of $\alpha_l = 2$ has been criticized since it often yields too smooth sample paths for many applications ([38]) and because computational difficulties can arise with this choice (see the Appendix). Thus $1 < \alpha_l < 2$ is typically chosen in the power exponential family ([4]), although the process is then not even once differentiable, sometimes not ideal for applications.

Another popular choice of the correlation function is the Matérn correlation. The isotropic, stationary form of the Handcock-Stein-Wallis parametrization of the Matérn function was introduced in [15, 16] and was extended to the non-stationary case in [27] via kernel convolution. When $\alpha_l = (2k+1)/2$ for $k \in \mathbb{N}$, the Matérn correlation has a closed-form expression. For example, when $\alpha_l = 1/2$, the Matérn correlation reduces to the power exponential correlation with $\alpha_l = 1$; when $\alpha_l \to \infty$, it reduces to Gaussian correlation. One nice feature of Matérn correlation is that its sample paths are $\lfloor \alpha_l - 1 \rfloor$ times differentiable, so the smoothness of the process can be directly controlled by the roughness parameters. Hence, it has become the recommended choice for the correlation function in spatial modeling ([38]). One of the most frequently used Matérn correlation functions is $\alpha_l = 5/2$, which has the form

$$c_l(d_l) = \left(1 + \frac{\sqrt{5}d_l}{\gamma_l} + \frac{5d_l^2}{3\gamma_l^2}\right)\exp\left(-\frac{\sqrt{5}d_l}{\gamma_l}\right), \tag{2.5}$$

where $d_l$ stands for any of the $|x_{il} - x_{jl}|$.

Use of Matérn correlation functions has been less popular in the computer model emulation literature. Here is an argument as to why (2.5) should be seriously considered for emulation, noting first that it is computationally tractable. Denoting $\tilde{d}_l = d_l/\gamma_l$, the following is easy to establish for (2.5).

- When $\tilde{d}_l \to 0$, $c_l(\tilde{d}_l) \approx 1 - C\tilde{d}_l^2$ with $C > 0$ being a constant. This thus behaves similarly to $\exp(-\tilde{d}_l^2) \approx 1 - \tilde{d}_l^2$, which corresponds to the power exponential correlation with $\alpha_l = 2$ (i.e., Gaussian correlation). This suggests that the Matérn correlation in (2.5) will maintain the smoothness induced by Gaussian correlation for nearby inputs.
- When $\tilde{d}_l \to \infty$, the dominant part of $c_l(\tilde{d}_l)$ is $\exp\left(-\sqrt{5}\tilde{d}_l\right)$ which matches the power exponential correlation with $\alpha_l = 1$. Thus the Matérn correlation in (2.5) prevents the correlation from decreasing quickly with distance, as does the Gaussian correlation. This can be of benefit in the computer model emulation since some inputs may have almost no effect on the computer model, which would correspond to near constant correlations even for distant inputs.

We have also found that the Matérn correlation function with $\alpha_l = 5/2$ yields very good empirical results in emulation. In addition, it is the default correlation function in the DiceKriging package. For these reasons, it will be used as the default correlation function for the numerical study in Section 5. However, our results are applicable to the much larger class of correlation functions listed in Table 1, as shown in Section 3.

### 2.2. *Marginal likelihood and marginal posterior.*

2.2.1. *Marginal likelihood.* Although maximum likelihood estimation of all parameters of the covariance function is possible, it has become standard to treat the mean parameters and variance in a fully objective Bayesian fashion, since they can be dealt with in closed-form in the Bayesian computations. Thus these parameters are assigned the objective prior

$$\pi(\boldsymbol{\theta}, \sigma^2) \propto \frac{1}{(\sigma^2)^a}\,.$$

with a fixed $a > 0$, where $a = 1$ corresponds to the standard reference prior. (It has become customary to also compare results with other choices of $a$, so we allow that in what follows.)

Using this prior to marginalize out the mean and variance parameters in the likelihood function, we obtain the marginal likelihood

$$L(\boldsymbol{\gamma} \mid \mathbf{y}^{\mathscr{D}}) \propto |\mathbf{R}|^{-\frac{1}{2}} |\mathbf{h}^T(\mathbf{x}^{\mathscr{D}})\mathbf{R}^{-1}\mathbf{h}(\mathbf{x}^{\mathscr{D}})|^{-\frac{1}{2}} \left(S^2\right)^{-(\frac{n-q}{2}+a-1)}, \tag{2.6}$$

where $\mathbf{h}(\mathbf{x}^{\mathscr{D}})$ is the $n \times q$ basis matrix with the $(i,j)$ entry $h_j(\mathbf{x}_i^{\mathscr{D}})$; $S^2 = (\mathbf{y}^{\mathscr{D}})^T\mathbf{Q}\mathbf{y}^{\mathscr{D}}$ with $\mathbf{Q} = \mathbf{R}^{-1}\mathbf{P}_{\boldsymbol{R}}$ and $\mathbf{P}_{\boldsymbol{R}} = \mathbf{I}_n - \mathbf{h}(\mathbf{x}^{\mathscr{D}})\{\mathbf{h}^T(\mathbf{x}^{\mathscr{D}})\mathbf{R}^{-1}\mathbf{h}(\mathbf{x}^{\mathscr{D}})\}^{-1}\mathbf{h}^T(\mathbf{x}^{\mathscr{D}})\mathbf{R}^{-1}$, with $\mathbf{I}_n$ being the identity matrix of size $n$.

Assuming the roughness parameters $\boldsymbol{\alpha} = (\alpha_1, \ldots, \alpha_p)$ have been prespecified, the range parameters of the correlation function can be estimated by maximizing (2.6), which is denoted as maximum marginal likelihood estimator (MMLE). While this approach was argued in [4] to be superior to maximum likelihood estimation (MLE), we will see that it is still non-robust, in the sense that will be defined in Section 3. The main problem is that the marginal likelihood will often not go to zero in the tails and, indeed, can be increasing. Thus it was argued in [23, 37] that the marginal likelihood needs to be augmented by the reference prior for the range parameters.

2.2.2. *Reference prior and posterior.* The reference prior for a separable product correlation function was developed in [28] and is given by

$$\pi^R(\boldsymbol{\theta}, \sigma^2, \boldsymbol{\gamma}) \propto \frac{\pi^R(\boldsymbol{\gamma})}{(\sigma^2)^a}\,, \tag{2.7}$$

with $\pi^R(\boldsymbol{\gamma}) \propto |\mathbf{I}^*(\boldsymbol{\gamma})|^{1/2}$, where $\mathbf{I}^*(\cdot)$ is the expected Fisher information matrix as below,

$$\mathbf{I}^*(\boldsymbol{\gamma}) = \begin{pmatrix} n-q & \mathrm{tr}(\mathbf{W}_1) & \mathrm{tr}(\mathbf{W}_2) & \ldots & \mathrm{tr}(\mathbf{W}_p) \\ & \mathrm{tr}(\mathbf{W}_1^2) & \mathrm{tr}(\mathbf{W}_1\mathbf{W}_2) & \ldots & \mathrm{tr}(\mathbf{W}_1\mathbf{W}_p) \\ & & \mathrm{tr}(\mathbf{W}_2^2) & \ldots & \mathrm{tr}(\mathbf{W}_2\mathbf{W}_p) \\ & & & \ddots & \vdots \\ & & & & \mathrm{tr}(\mathbf{W}_p^2) \end{pmatrix}, \tag{2.8}$$

where $\mathbf{W}_l = \dot{\mathbf{R}}_l\mathbf{Q}$, for $1 \le l \le p$, and $\dot{\mathbf{R}}_l$ is the partial derivative of the correlation matrix $\mathbf{R}$ with respect to the $l^{th}$ range parameter.

The marginal posterior of $\boldsymbol{\gamma}$ with regard to this reference prior is

$$p(\boldsymbol{\gamma} \mid \mathbf{y}^{\mathscr{D}}) \propto L(\boldsymbol{\gamma} \mid \mathbf{y}^{\mathscr{D}})\, |\mathbf{I}^*(\boldsymbol{\gamma})|^{1/2}. \tag{2.9}$$

Sampling from this posterior requires a Metropolis-type algorithm and each evaluation of the likelihood typically requires $O(n^3)$ flops for the inverse of the correlation matrix, which is computationally prohibitive for many applications. Moreover, the computation error can be very large when the correlation matrix is close to the matrix of all ones. For these reasons, it is common ([4, 37]) to instead simply estimate $\boldsymbol{\gamma}$ by its marginal posterior mode, using (2.9),

$$(\hat{\gamma}_1, \ldots \hat{\gamma}_p) = \underset{\gamma_1, \ldots, \gamma_p}{\mathrm{argmax}} \left\{ L(\gamma_1, \ldots, \gamma_p \mid \mathbf{y}^{\mathscr{D}})\, \pi^R(\gamma_1, \ldots, \gamma_p) \right\}. \tag{2.10}$$

2.2.3. *Parameterizations.* Maximum likelihood estimation is invariant under the choice of parameterization, but the posterior mode is not invariant because of the presence of the Jacobian for the prior. Here are three common ways of parameterizing the range parameters in the power exponential correlation function ([5, 28, 3, 4, 37]), for any $l = 1, \cdots, p$:

$$c_{\gamma_l}(|x_{il} - x_{jl}|) = \exp\{-(|x_{il} - x_{jl}|/\gamma_l)^{\alpha_l}\}, \tag{2.11}$$

$$c_{\tilde{\beta}_l}(|x_{il} - x_{jl}|) = \exp\{-\tilde{\beta}_l |x_{il} - x_{jl}|^{\alpha_l}\}, \tag{2.12}$$

$$c_{\tilde{\xi}_l}(|x_{il} - x_{jl}|) = \exp\left\{-\exp(\tilde{\xi}_l)|x_{il} - x_{jl}|^{\alpha_l}\right\}. \tag{2.13}$$

Table 1 gives various correlation functions in their natural parameterizations, in which the range parameter and roughness parameter are independent; we will call this the $\alpha$-*free* parameterization of the range parameter. In contrast, in the above parameterizations of the power exponential correlation function, $\tilde{\beta}_l = \gamma_l^{-\alpha_l}$ and $\tilde{\xi}_l = \log(\gamma_l^{-\alpha_l})$ both depend on $\alpha_l$. We will also consider the following transformations of the $\alpha$-*free* parameterization (dropping the subscript $l$ for convenience).

Definition 2.1. *For the range parameters $\gamma$ in Table 1,*
*(i) $\beta = 1/\gamma$ will be called the inverse range parameter;*
*(ii) $\xi = \log(1/\gamma)$ will be called the log inverse range parameter.*

Note that $\tilde{\beta} = \beta^\alpha$ and $\tilde{\xi} = \alpha\xi$. The mode of the posterior distributions for the $\tilde{\xi}$ and $\xi$ parameterizations will be the same (properly transformed), because the Jacobians of the transformations differ only by the prefixed constant $\alpha$; thus we need to consider only the $\xi$ — and not the $\tilde{\xi}$ — parameterization of the power exponential correlation function in what follows. On the other hand, the posterior modes of $\tilde{\beta}$ and $\beta$ are not the same (when transformed), so we have to consider both parameterizations in what follows.

2.2.4. *Predictions using the emulator.* After obtaining the estimates of the range parameters under a specified parameterization, transform back to obtain the corresponding $\boldsymbol{\gamma} = (\gamma_1, \ldots, \gamma_p)$, after which the predictive distribution of $y(\mathbf{x}^*)$, given $\mathbf{y}^{\mathscr{D}}$ and $\boldsymbol{\gamma}$, is a Student's $t$-distribution,

$$y(\mathbf{x}^*) \mid \mathbf{y}^{\mathscr{D}}, \boldsymbol{\gamma} \sim t(\hat{y}(\mathbf{x}^*), \hat{\sigma}^2 c^{**}, n-q)\,, \tag{2.14}$$

with $n-q$ degrees of freedom, where

$$\begin{aligned}\hat{y}(\mathbf{x}^*) &= \mathbf{h}(\mathbf{x}^*)\hat{\boldsymbol{\theta}} + \mathbf{r}^T(\mathbf{x}^*)\mathbf{R}^{-1}\left(\mathbf{y}^{\mathscr{D}} - \mathbf{h}(\mathbf{x}^{\mathscr{D}})\hat{\boldsymbol{\theta}}\right),\\ \hat{\sigma}^2 &= (n-q)^{-1}\left(\mathbf{y}^{\mathscr{D}} - \mathbf{h}(\mathbf{x}^{\mathscr{D}})\hat{\boldsymbol{\theta}}\right)^T\mathbf{R}^{-1}\left(\mathbf{y}^{\mathscr{D}} - \mathbf{h}(\mathbf{x}^{\mathscr{D}})\hat{\boldsymbol{\theta}}\right),\\ c^{**} &= c(\mathbf{x}^*,\mathbf{x}^*) - \mathbf{r}^T(\mathbf{x}^*)\mathbf{R}^{-1}\mathbf{r}(\mathbf{x}^*) + \left(\mathbf{h}(\mathbf{x}^*) - \mathbf{h}^T(\mathbf{x}^{\mathscr{D}})\mathbf{R}^{-1}\mathbf{r}(\mathbf{x}^*)\right)^T\\ &\quad \times\left(\mathbf{h}^T(\mathbf{x}^{\mathscr{D}})\mathbf{R}^{-1}\mathbf{h}(\mathbf{x}^{\mathscr{D}})\right)^{-1}\left(\mathbf{h}(\mathbf{x}^*) - \mathbf{h}^T(\mathbf{x}^{\mathscr{D}})\mathbf{R}^{-1}\mathbf{r}(\mathbf{x}^*)\right),\end{aligned}$$

with $\hat{\boldsymbol{\theta}} = \left(\mathbf{h}^T(\mathbf{x}^{\mathscr{D}})\mathbf{R}^{-1}\,\mathbf{h}(\mathbf{x}^{\mathscr{D}})\right)^{-1}\mathbf{h}^T(\mathbf{x}^{\mathscr{D}})\mathbf{R}^{-1}\mathbf{y}^{\mathscr{D}}$ being the generalized least squares estimator for $\boldsymbol{\theta}$; $\mathbf{R}$ being the correlation matrix corresponding to the design inputs and $\mathbf{r}(\mathbf{x}^*) = \left(c(\mathbf{x}^*,\mathbf{x}_1^{\mathscr{D}}),\cdots,c(\mathbf{x}^*,\mathbf{x}_n^{\mathscr{D}})\right)^T$, both obtained by plugging in the estimated $\boldsymbol{\gamma}$ values. The corresponding prediction and any quantile of the predictive distribution are then readily available.

2.3. *Profile likelihood.* For comparison purposes, we will also consider the full likelihood approach, which utilizes the MLE for the mean and variance parameters, $\hat{\boldsymbol{\theta}}_{MLE} = \hat{\boldsymbol{\theta}}$, $\hat{\sigma}^2_{MLE} = (n-q)\hat{\sigma}^2/n$, where $\hat{\boldsymbol{\theta}}$ and $\hat{\sigma}^2$ are defined in (2.14). Plugging $\hat{\boldsymbol{\theta}}_{MLE}$ and $\hat{\sigma}^2_{MLE}$ into (2.2) and ignoring the normalizing constant, the likelihood of (2.2) reduces to the profile likelihood

$$L(\boldsymbol{\gamma}\mid\mathbf{y}^{\mathscr{D}},\hat{\sigma}^2_{MLE},\hat{\boldsymbol{\theta}}_{MLE}) \propto |\mathbf{R}|^{-\frac{1}{2}}(S^2)^{-\frac{n}{2}}\,. \tag{2.15}$$

To complete the MLE analysis, $\boldsymbol{\gamma}$ is estimated by the mode of this profile likelihood and denoted by $\hat{\boldsymbol{\gamma}}_{MLE}$. The predictive distribution of a new input $\mathbf{x}^*$, conditional on the previous outputs and the MLE, is

$$y(\mathbf{x}^*)\mid\mathbf{y}^{\mathscr{D}},\hat{\sigma}^2_{MLE},\hat{\boldsymbol{\theta}}_{MLE},\hat{\boldsymbol{\gamma}}_{MLE} \sim N(\hat{y}_{MLE}(\mathbf{x}^*),\hat{\sigma}^2_{MLE}c^*_{MLE}), \tag{2.16}$$

where $\hat{y}_{MLE}(\mathbf{x}^*) = \hat{y}(\mathbf{x}^*)$, with $\hat{y}(\mathbf{x}^*)$ defined in (2.14), and $c^*_{MLE} = c_{MLE}(\mathbf{x}^*,\mathbf{x}^*) - \mathbf{r}^T_{MLE}(\mathbf{x}^*)\mathbf{R}^{-1}_{MLE}\mathbf{r}_{MLE}(\mathbf{x}^*)$, obtained by plugging $\hat{\boldsymbol{\gamma}}_{MLE}$ into $c_{MLE}(\mathbf{x}^*,\mathbf{x}^*)$, $\mathbf{r}_{MLE}(\mathbf{x}^*)$ and $\mathbf{R}_{MLE}$.

The profile likelihood is sometimes very flat in the tails, resulting in $\hat{\boldsymbol{\gamma}}_{MLE}$ being near zero and $\hat{\mathbf{R}}_{MLE}$ being near $\mathbf{I}_n$ (see the details in Section 3). This can be shown to result in the predicted mean, $\hat{y}_{MLE}(\mathbf{x}^*)$, being essentially an impulse function at each of the observations, while following the GaSP mean elsewhere. Figure 1 gives an example of this scenario, where the GaSP mean is assumed to be a constant. In the left panel of the figure, the roughness parameter was $\alpha = 1$ for the power exponential correlation function, and both the MLE and MMLE became essentially degenerate, while the prediction from the posterior mode approach was reasonable (although not quite smooth enough). In the right panel of the figure, the roughness parameter was $\alpha = 1.9$; here both the MMLE and marginal posterior mode approaches gave excellent predictions, but the profile likelihood approach still resulted in a degenerate prediction. Such degeneracies are somewhat unusual in one-dimension, but are not particularly unusual with higher dimensional inputs, as shown numerically in Section 5.

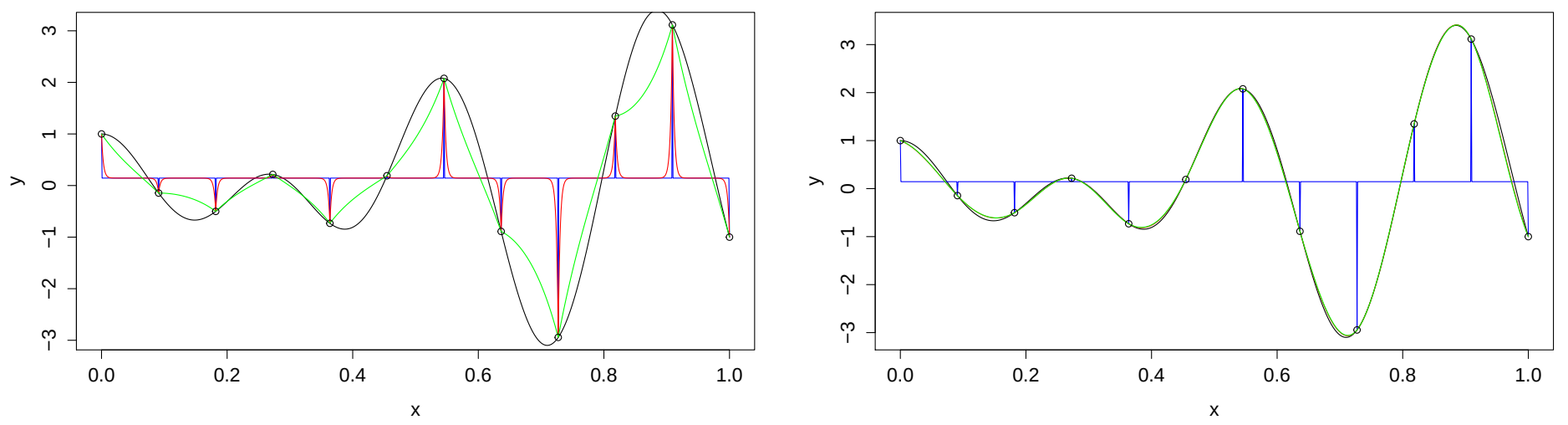

FIG 1. *Emulation of the function* $y = 3sin(5\pi x)x + cos(7\pi x)$*, graphed as the black solid curves (overlapping the green and red curves in the right panel). The design for the input* $x$ *is equally spaced from* $[0, 1]$ *with* $n = 12$*, with the resulting function values indicated by the black circles. A constant mean function is used. The left panel is for* $\alpha = 1$ *and the right panel for* $\alpha = 1.9$*, for the power exponential correlation function. The blue curves (which are essentially unit impulse functions at the observations and constant elsewhere) give the emulator mean obtained from the profile likelihood approach; the red curves give the emulator mean from the MMLE approach; and the green curves give the emulator mean arising from the maximum posterior mode approach with the reference prior.*

**3. Robust parameter estimation for GaSP Models.** In this section, we explore the ways in which GaSP emulator construction can fail, developing the "robustness criteria" that are needed to avoid such failures. We then examine which estimation methods satisfy the criteria. To begin, it is pedagogically useful to look at a special case ([23]), where the analysis is essentially closed-form. The proofs of the lemmas and theorems in this section are provided in the supplementary materials ([14]).

3.1. *A closed-form example for the profile likelihood and marginal likelihood.* Suppose the input is one-dimensional and that the design is equally spaced with the design points being $d_0$ units apart. Consider a constant mean $h(x) = 1$ and power exponential correlation with roughness parameter $\alpha = 1$. Denote $\rho = e^{-d_0/\gamma}$, write $c(x_i, x_j) = \rho^{\Delta_{ij}}$, with $\Delta_{ij} = |x_i - x_j|/d_0$, and write $y(x_i^{\mathscr{D}})$ as $y_i$ to simplify the notation. The closed-form logarithm of the profile likelihood and marginal likelihood (obtained by integrating out the mean and variance parameters using the standard reference prior), as well as their limiting values when $\rho \to 0$ and $\rho \to 1$, are given in the supplementary materials ([14]). From these, we can establish the following condition, under which the mode of the profile likelihood occurs at $\rho = 0$.

LEMMA 3.1. *A necessary and sufficient condition that the mode of the profile likelihood in (2.15) is at* $\rho = 0$ [causing the unwelcome degeneracy] *is, defining* $\bar{y} = \sum_{i=1}^{n} y_i/n$,

$$\sum_{i=1}^{n-1}(y_i - \bar{y})(y_{i+1} - \bar{y}) \leq 0\,. \tag{3.1}$$

The intuition behind Lemma 3.1 comes from the fact that, in this case, the GaSP becomes an autoregressive model of order 1. When the empirical lag-1 autocorrelation is less than zero, the profile likelihood estimate of the correlation $\rho$ will be zero, since the correlation $\rho$ is parameterized to be nonnegative here. On the other hand, if either likelihood is maximized at $\rho = 1$, then $\mathbf{R} = \mathbf{1}_n\mathbf{1}_n^T$, where $\mathbf{1}_n$ is the vector of all ones, so that the correlation matrix becomes ill-conditioned, causing large approximation errors in computation of its inverse.

For the general case considered in the remainder of the paper, explicit results such as that in Lemma 3.1 are not available. However, we can still look at the tail rates (corresponding to $\rho$ going to 0 or 1) for various likelihoods and posteriors and assess when problems will occur. We formalize these notions in the next subsection, through our criteria for robust estimation.

3.2. *Robust estimation.* As discussed in the previous section, when $\mathbf{R} \approx \mathbf{I}_n$, the GaSP predictive mean will degenerate to the fitted mean and impulse functions at the observed inputs, as happened in Figure 1. When $\mathbf{R} \approx \mathbf{1}_n\mathbf{1}_n^T$, the correlation matrix $\mathbf{R}$ is almost singular, leading to very large computational errors in the GaSP predictive mean. Robust estimation of the parameters is defined as avoiding these two possible problems.

DEFINITION 3.1. *(Robust Estimation.) Estimation of the parameters in the GaSP is called robust, if the following two situations do NOT happen:*

*(i)* $\hat{\mathbf{R}} = \mathbf{1}_n\mathbf{1}_n^T$,
*(ii)* $\hat{\mathbf{R}} = \mathbf{I}_n$,

*where* $\hat{\mathbf{R}}$ *is the estimated correlation matrix.*

Note that the predictive mean of the GaSP is not well-defined in these two situations when the inputs are at one of the design points, but it can be defined as the limit as $\hat{\mathbf{R}} \to \mathbf{1}_n \mathbf{1}_n^T$, and $\hat{\mathbf{R}} \to \mathbf{I}_n$.

The following basic lemma is immediate from the definition of the correlation matrix.

LEMMA 3.2. *Robustness is lacking in either of the following two cases.*
*Case 1. If, for all* $1 \leq l \leq p$, $\hat{\gamma}_l = \infty$ *(or* $\hat{\xi}_l = -\infty$ *or* $\hat{\beta}_l = 0$ *in the other parameterizations), then* $\hat{\mathbf{R}} = \mathbf{1}_n \mathbf{1}_n^T$.
*Case 2. If* $\exists l$, $1 \leq l \leq p$, *for which* $\hat{\gamma}_l = 0$ *(equivalent to* $\hat{\xi}_l = \infty$ *or* $\hat{\beta}_l = \infty$*), then* $\hat{\mathbf{R}} = \mathbf{I}_n$.

Note that it is generally fine if some (but not all) of the estimated $\gamma_l$ are close to $\infty$, because this will just make $\hat{\mathbf{R}}_l \approx \mathbf{1}_n \mathbf{1}_n^T$ for some $l$ but not $\hat{\mathbf{R}} \approx \mathbf{1}_n \mathbf{1}_n^T$. In such a situation, the inputs associated with the large $\gamma_l$ can be called *inert* inputs, since they will have only a small effect on the outputs. Indeed, this is a desirable situation, since such inputs could be removed from the emulator, simplifying and improving the approximation.

The MLE, MMLE and marginal posterior modes (for the various parameterizations) all reduce to mode estimation with regard to a function $G(\boldsymbol{\gamma})$. Thus the following guarantees that the problematic situations cannot occur.

COROLLARY 3.1. *Estimation of* $\boldsymbol{\gamma} = (\gamma_1, \ldots, \gamma_p)^T$ *as the mode of a non-negative function* $G(\boldsymbol{\gamma})$ *is robust if* $G(\boldsymbol{\gamma}) \to 0$, *under the following two situations:*

*(i)* $\exists\ l$, $1 \leq l \leq p$, $\gamma_l \to 0$,
*(ii) For all* $l$, $1 \leq l \leq p$, $\gamma_l \to \infty$.

COROLLARY 3.2. *Estimation of any monotonic transformation of the range parameters* $\boldsymbol{\zeta} = \mathbf{f}(\boldsymbol{\gamma}) = (f_1(\boldsymbol{\gamma}), \ldots, f_p(\boldsymbol{\gamma}))^T$, *by the mode of its marginal posterior, is robust if*

$$L(\mathbf{f}^{-1}(\boldsymbol{\zeta}) \mid \mathbf{y}^{\mathscr{D}}) \pi^R(\mathbf{f}^{-1}(\boldsymbol{\zeta})) \left| \frac{\partial \mathbf{f}^{-1}(\boldsymbol{\zeta})}{\partial \boldsymbol{\zeta}} \right| \to 0$$

*under the following two situations:*

*(i)* $\exists\ l$, $1 \leq l \leq p$, $f_l^{-1}(\boldsymbol{\zeta}) \to 0$,
*(ii) For all* $l$, $1 \leq l \leq p$, $f_l^{-1}(\boldsymbol{\zeta}) \to \infty$.

*where* $\mathbf{f}^{-1}(\boldsymbol{\zeta}) = (f_1^{-1}(\boldsymbol{\zeta}), ..., f_p^{-1}(\boldsymbol{\zeta}))^T$.

3.3. *Robustness Results.* From the results in the previous section, it is clear that we should compute the tail rates, in terms of $\boldsymbol{\gamma}$, of the marginal likelihood, profile likelihood, and the various posteriors to see if they are robust. Computation of the tail rates of the posteriors requires computation of the tail rates of the reference prior, as well as the tail rates of the marginal likelihood. We need the following two mild assumptions (c.f., [5, 32]) to establish the main results concerning these rates.

ASSUMPTION 3.1. *For any* $d_l \geq 0$ *and* $1 \leq l \leq p$, $c_l(d_l) = c_l^0(d_l/\gamma_l)$, *where* $c_l^0(\cdot)$ *is a correlation function that satisfies* $\lim\limits_{u\to\infty} c_l^0(u) = 0$.

ASSUMPTION 3.2. *For any* $1 \leq l \leq p$, *as* $\gamma_l \to \infty$,

$$\mathbf{R}_l(\gamma_l) = \mathbf{1}_n \mathbf{1}_n^T + \nu_l(\gamma_l)\mathbf{D}_l + \nu_l(\gamma_l)\omega_l(\gamma_l)(\mathbf{D}_l^* + \mathbf{B}_l(\gamma_l)),$$

*where* $\mathbf{D}_l$ *is a nonsingular and symmetric matrix with* $\mathbf{1}_n^T\mathbf{D}_l^{-1}\mathbf{1}_n \neq 0$, $\mathbf{D}_l^*$ *is a fixed matrix,* $\nu_l(\gamma_l) > 0$ *is a non-increasing and differentiable function,* $\omega_l(\gamma_l)$ *is a differentiable function, and* $\mathbf{B}_l(\gamma_l)$ *is a differentiable matrix (incorporating the higher order terms of the expansion), satisfying*

$$\nu_l(\gamma_l) \to 0,\, \omega_l(\gamma_l) \to 0,\, \frac{\omega_l'(\gamma_l)}{\frac{\partial}{\partial\gamma_l}\log\, \nu_l(\gamma_l)} \to 0,\, ||\mathbf{B}_l(\gamma_l)||_\infty \to 0,\, \frac{||\frac{\partial}{\partial\gamma_l}\mathbf{B}_l(\gamma_l)||_\infty}{\frac{\partial}{\partial\gamma_l}\log(\omega_l(\gamma_l))} \to 0,$$

*where* $\omega_l'(\gamma_l) = \partial\omega_l(\gamma_l)/\partial\gamma_l$, *and* $\|\boldsymbol{B}\|_\infty = \max_{i,j}|a_{ij}|$ *with* $a_{ij}$ *being the* $(i,j)$ *entry of the matrix* $\boldsymbol{B}$.

The first assumption ensures that the correlation function will decrease to zero as the distance between two points goes to infinity. The second assumption guarantees that the first two small terms in the Taylor expansion of the correlation function decrease to zero at an appropriate rate as $\gamma_l \to \infty$. The assumptions hold for all the correlation functions listed in Table 1, in which the functions $\nu_l$ and $\omega_l$ are also given.

The following lemma gives the tail rates for the marginal and profile likelihoods.

LEMMA 3.3. *(Tail rates of the marginal likelihood and profile likelihood.) If Assumption 3.1 and Assumption 3.2 hold for each of the* $\mathbf{R}_l$, $1 \leq l \leq p$, *the marginal likelihood and profile likelihood have the following tail rates.*

*(i) If* $\exists l$, $1 \leq l \leq p$, *such that* $\gamma_l \to 0$, *the marginal likelihood and profile likelihood both exist and are greater than zero.*

*(ii) If $\gamma_l \to \infty$ for all $l$, $1 \leq l \leq p$, and $\mathcal{C}(\mathbf{h}(\mathbf{x}^{\mathscr{D}}))$ denotes the column space of the mean basis matrix $\mathbf{h}(\mathbf{x}^{\mathscr{D}})$, the marginal likelihood satisfies*

$$L(\boldsymbol{\gamma} \mid \mathbf{y}^{\mathscr{D}}) = \begin{cases} O\left( (\sum_{l=1}^{p} \nu_l(\gamma_l))^{a-1/2} \right), & \mathbf{1}_n \notin \mathcal{C}(\mathbf{h}(\mathbf{x}^{\mathscr{D}})), \\ O\left( (\sum_{l=1}^{p} \nu_l(\gamma_l))^{a-1} \right), & \mathbf{1}_n \in \mathcal{C}(\mathbf{h}(\mathbf{x}^{\mathscr{D}})) \,. \end{cases}$$

*The profile likelihood, in this case, satisfies*

$$L(\boldsymbol{\gamma} \mid \mathbf{y}^{\mathscr{D}}, \hat{\sigma}^2_{MLE}, \hat{\boldsymbol{\theta}}_{MLE}) = O\left( (\sum_{l=1}^{p} \nu_l(\gamma_l))^{1/2} \right).$$

Part (i) of this lemma indicates that the marginal likelihood and profile likelihood could have their modes at $\mathbf{R} = \mathbf{I}_n$ and thus could potentially be non-robust; one such case was given in Figure 1.

Part (ii) of the lemma shows that the mode of the marginal likelihood could be at $\mathbf{R} = \mathbf{1}_n \mathbf{1}_n^T$ for the frequently used setting of $a = 1$ and $\mathbf{1}_n \in \mathcal{C}(\mathbf{h}(\mathbf{x}^{\mathscr{D}}))$. On the other hand, the profile likelihood will decrease to zero at this limit, so it cannot be non-robust in this fashion. A byproduct of Lemma 3.3 is that, when $a = 1$ and $\mathbf{1}_n \in \mathcal{C}(\mathbf{h}(\mathbf{x}^{\mathscr{D}}))$, use of a constant prior for $\boldsymbol{\gamma}$ would result in an improper posterior distribution, consistent with the result for isotropic case given in [5].

The asymptotic behaviors of the reference prior for the two limiting cases of interest are given in the lemma below.

LEMMA 3.4. *(Tail rates of the prior.) If Assumption 3.1 and Assumption 3.2 hold for each of the $\mathbf{R}_l$, $1 \leq l \leq p$, then $\pi^R(\boldsymbol{\gamma})$ has the following two limiting properties. Here $\boldsymbol{\gamma}_E$ denotes the vector of $\gamma_l$ for all $l \in E$, $E \subset \{1, 2, ..., p\}$, and $\boldsymbol{\gamma}_{-E}$ denotes the complementary vector.*

*(i) As $\boldsymbol{\gamma}_E \to \mathbf{0}$,*

$$\pi^R(\boldsymbol{\gamma}) \leq C(\boldsymbol{\gamma}_{-E}) \left[ \prod_{l \in E} \operatorname{tr} \left( \frac{\partial \mathbf{R}}{\partial \gamma_l} \right)^2 \right]^{1/2},$$

*where $C(\boldsymbol{\gamma}_{-E})$ is constant in $\boldsymbol{\gamma}_E$.*

*(ii) As $\gamma_l \to \infty$ for all $l$, $1 \leq l \leq p$, if $\mathbf{1} \notin \mathcal{C}(\boldsymbol{h}(\boldsymbol{x}))$,*

$$\pi^R(\boldsymbol{\gamma}) \leq C_1 \left| \frac{\prod_{l=1}^{p} \nu_l'(\gamma_l)}{(\sum_{l=1}^{p} \nu_l(\gamma_l))^p} \right|, \tag{3.2}$$

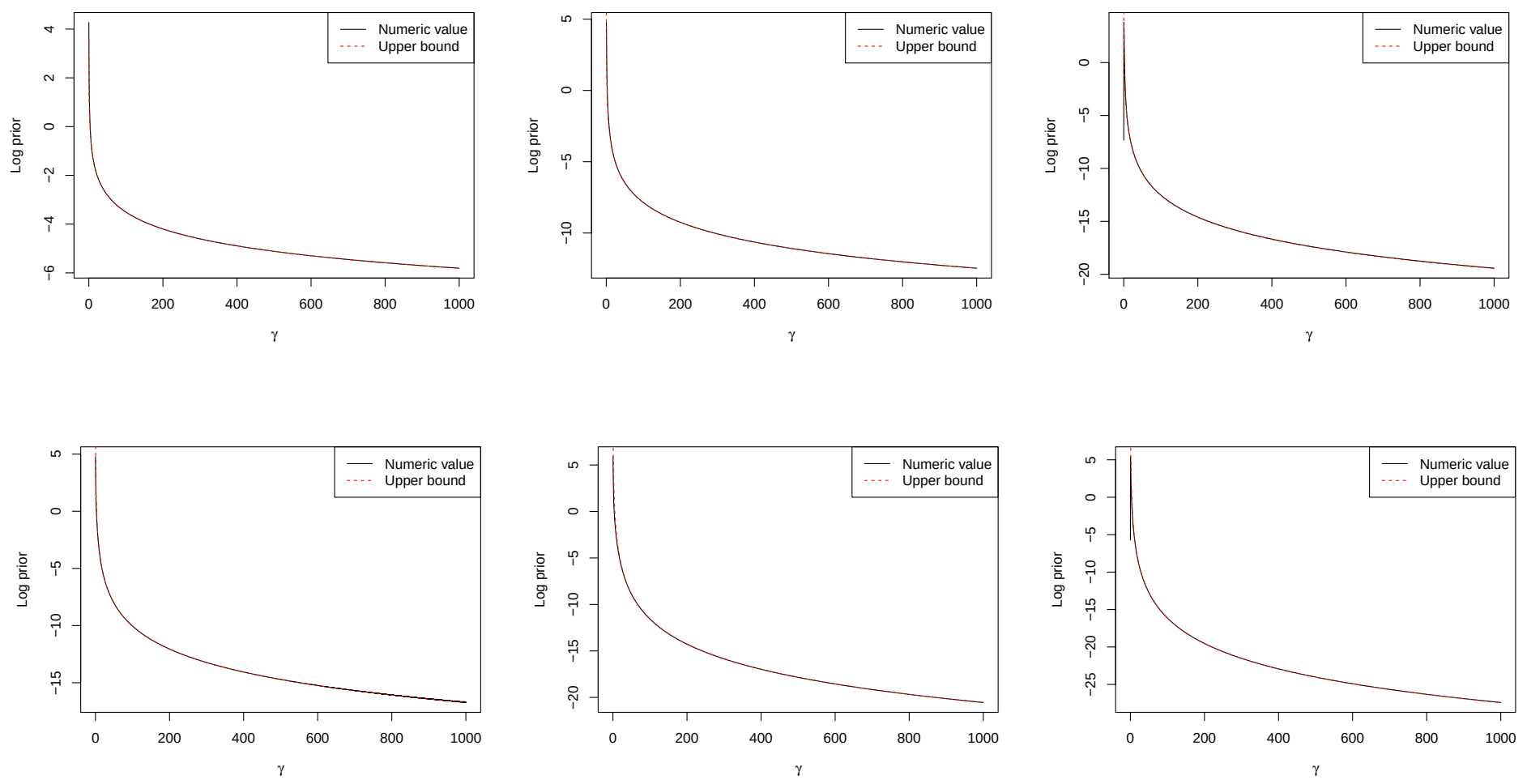

Fig 2. *The tail behavior of the reference prior (black curves), and its upper bound (red curves) from Lemma 3.4 part (ii), when $\gamma_1 = \ldots = \gamma_p \to \infty$. The power exponential correlation function is used with fixed $\alpha_l = 1.9$, $1 \le l \le p$. The first row is for the case in which $\mathbf{1}_n \notin \mathcal{C}(\mathbf{h}(\mathbf{x}^{\mathscr{D}}))$, while the second row is for $\mathbf{1}_n \in \mathcal{C}(\mathbf{h}(\mathbf{x}^{\mathscr{D}}))$. From left to right, the dimension of the inputs are $p = 1$, $p = 2$ and $p = 3$. The prior and bounds are evaluated at points uniformly sampled from $[0, 1]^p$. The black curves and red curves overlap when $\gamma_l$ is large.*

*where $\nu_l'(\gamma_l) = \partial \nu_l(\gamma_l)/\partial \gamma_l$; if $\mathbf{1} \in \mathcal{C}(\mathbf{h}(\mathbf{x}^{\mathscr{D}}))$ and $p \ge 2$,*

$$\pi^R(\boldsymbol{\gamma}) \le C_2 \left| \frac{\prod_{l=1}^{p} \nu_l'(\gamma_l)}{(\sum_{l=1}^{p} \nu_l(\gamma_l))^p} \right| \left| \sum_{l=1}^{p} \frac{\nu_l^2(\gamma_l)\omega_l'(\gamma_l)}{\nu_l'(\gamma_l)\nu_m(\gamma_m)} \right|,$$

*for every index $m$ between 1 and $p$; if $\mathbf{1} \in \mathcal{C}(\mathbf{h}(\mathbf{x}^{\mathscr{D}}))$ and $p = 1$,*

$$\pi^R(\boldsymbol{\gamma}) \le C_3 |\omega_1'(\gamma_1)|.$$

*where $C_1$, $C_2$ and $C_3$ are all positive and not related to $\gamma_l$.*

The bounds for the one-dimensional case in Lemma 3.4 (ii) were proved in [5]. These results are a generalization of the $p$ dimensional results in [28], which considered only separable designs.

Interestingly, the bounds in part (ii) of Lemma 3.4 seem to be almost exact in numerical examples we have studied for the power exponential correlation function. Figure 2 presents some of the evidence for this.

The following theorem states that, under the $\boldsymbol{\gamma}$ and $\boldsymbol{\xi}$ parameterizations and when $a = 1$, the mode of the marginal posterior with the reference prior for the range parameters will typically be robust for the correlation functions listed in Table 1. Similar theorems can be stated for other choices of $a$ but, since $a = 1$ is the near universal choice, we restrict the statement of the results to that case.

THEOREM 3.1. *Under the parameterizations of the range parameter* $\boldsymbol{\gamma}$ *and log inverse range* $\boldsymbol{\xi}$ *in Definition 2.1, the posterior mode in (2.9) with* $a = 1$ *is robust for the product form of the power exponential, spherical, and Matérn correlation functions over the domain of* $\boldsymbol{\alpha}$ *listed in Table 1. In addition, the posterior mode of* $\boldsymbol{\gamma}$ *is robust for the rational quadratic correlation if* $\alpha_l > 1/2$, $1 \leq l \leq p$ *and the posterior mode of* $\boldsymbol{\xi}$ *is robust for the rational quadratic correlation over the entire domain of* $\boldsymbol{\alpha}$*.*

PROOF. Theorem 3.1 can be proved by verifying Corollary 3.1 and Corollary 3.2 using the results from Lemma 3.3 and Lemma 3.4. □

While use of the mode of the marginal posterior for the $\boldsymbol{\gamma}$ and $\boldsymbol{\xi}$ parameterizations is robust, the mode of the marginal posterior under other parameterizations, such as the $\tilde{\boldsymbol{\beta}}$ parameterization in (2.12), can be non-robust. Indeed, directly applying Lemma 3.4 and Lemma 3.3, the bounds on the tail rates of the marginal posterior under the various parameterizations (and also for the profile and marginal likelihood) are given in Table 2. For simplicity, we assume roughness parameters are kept the same, i.e. $\alpha_1 = \alpha_2 = \cdots = \alpha_p = \alpha$. The blue highlighted entries are those in which the tail rate is constant, so that there is a potential problem of non-robustness.

The red highlighted entries in Table 2 are quite surprising, as here the marginal posterior density becomes infinite in the tail, so that the mode will be at the problematical $\mathbf{1}_n\mathbf{1}_n^T$. (The following Corollary 3.3 establishes that there is no other infinite mode.) That the posterior mode for the $\tilde{\boldsymbol{\beta}}$ parameterization has this bizarre behavior has not been previously recognized, and should clearly rule out use of this parameterization (at least when estimating by the marginal posterior mode with the standard reference prior). Figure 3 gives numerical evidence of this feature, where we plot the log-marginal posterior as a function of $\tilde{\beta}_1 = \tilde{\beta}_2 = \tilde{\beta}$. Both examples have local modes with a finite marginal posterior, while the real modes with infinite posterior density occur as $\tilde{\beta}_1 = \tilde{\beta}_2 \to 0$.

The following lemma is needed to establish posterior propriety in the next subsection and also to establish Corollary 3.3. It calculates the tail rates when some, but not all, of the range parameters are close to zero.

| | $\mathbf{1}_n \in \mathcal{C}(\mathbf{h}(\mathbf{x}^{\mathscr{D}}))$ | | $\mathbf{1}_n \notin \mathcal{C}(\mathbf{h}(\mathbf{x}^{\mathscr{D}}))$ | |
|---|---|---|---|---|
| | $l \in E$, $\gamma_l \to 0$ | $\gamma_l \to \infty$ for all $l$ | $l \in E$, $\gamma_l \to 0$ | $\gamma_l \to \infty$ for all $l$ |
| Profile Lik | $O(1)$ | $O(\gamma_{(1)}^{-\alpha/2})$ | $O(1)$ | $O(\gamma_{(1)}^{-\alpha/2})$ |
| Marginal Lik | $O(1)$ | $O(1)$ | $O(1)$ | $O(\gamma_{(1)}^{-\alpha/2})$ |
| Post $\boldsymbol{\gamma}$, $p=1$ | $O(\frac{\exp(-C/\gamma^{\alpha})}{\gamma^{(\alpha+1)}})$ | $O(\gamma^{-\alpha-1})$ | $O(\frac{\exp(-C/\gamma^{\alpha})}{\gamma^{(\alpha+1)}})$ | $O(\gamma^{-\alpha/2-1})$ |
| $p \geq 2$ | $O(\prod_{l\in E} \frac{\exp(-C_l/\gamma_l^{\alpha})}{\gamma_l^{(\alpha+1)}})$ | $O(\frac{\prod_{l=1}^{p} \gamma_l^{-\alpha-1}}{\gamma_{(1)}^{(1-p)\alpha}})$ | $O(\prod_{l\in E} \frac{\exp(-C_l/\gamma_l^{\alpha})}{\gamma_l^{(\alpha+1)}})$ | $O(\frac{\prod_{l=1}^{p} \gamma_l^{-\alpha-1}}{\gamma_{(1)}^{(1/2-p)\alpha}})$ |
| Post $\tilde{\boldsymbol{\beta}}$, $p=1$ | $O(\exp(-\tilde{\beta}C))$ | $O(1)$ | $O(\exp(-\tilde{\beta}C))$ | $O(\tilde{\beta}^{-1/2})$ |
| $p \geq 2$ | $O(\prod_{l\in E} \exp(-\tilde{\beta}_l C_l))$ | $O(\tilde{\beta}_{(p)}^{-(p-1)})$ | $O(\prod_{l=1}^{p} \exp(-\tilde{\beta}_l)C_l)$ | $O(\tilde{\beta}_{(p)}^{-(p-1/2)})$ |
| Post $\tilde{\boldsymbol{\xi}}$, $p=1$ | $O(\exp(-\exp(\tilde{\xi})C+\tilde{\xi}))$ | $O(\exp(\tilde{\xi}))$ | $O(\exp(-\exp(\tilde{\xi})C))$ | $O(\exp(\tilde{\xi}/2))$ |
| $p \geq 2$ | $O(\prod_{l\in E} \exp(-\exp(\tilde{\xi}_l)C_l+\tilde{\xi}_l))$ | $O(\frac{\exp(\sum_{l=1}^{p-1} \tilde{\xi}_l)}{\exp((p-2)\tilde{\xi}_{(p)})})$ | $O(\prod_{l\in E} \exp(-\exp(\tilde{\xi}_l)C_l))$ | $O(\frac{\exp(\sum_{i=1}^{p-1} \tilde{\xi}_l)}{\exp((p-1/2)\tilde{\xi}_{(p)})})$ |

TABLE 2

*Tail behaviors of the profile likelihood, the marginal likelihood and the posterior distributions for different parameterizations of the power exponential correlation function, using the reference prior in (2.7) with $a = 1$. In the 2nd and 4th columns, $E$ is a nonempty set such that for $l \in E$, $\gamma_l \to 0$ (equivalent to $\tilde{\beta}_l \to \infty$ or $\tilde{\xi}_l \to \infty$), and $C$ and $C_l$ are positive numbers depending on $|x_{il}^{\mathscr{D}} - x_{jl}^{\mathscr{D}}|$, $1 \leq i, j \leq n$, $l \in E$. In the 3rd and 5th columns, $\gamma_l \to \infty$ (equivalent to $\tilde{\beta}_l \to 0$ or $\tilde{\xi}_l \to -\infty$), for all $1 \leq l \leq p$; in the stated tail rates, $\gamma_{(1)}$ is defined as the minimum of the $\gamma_l$, $\tilde{\beta}_{(p)}$ is the largest $\tilde{\beta}_l$, and $\tilde{\xi}_{(p)}$ is the largest $\tilde{\xi}_l$, where $1 \leq l \leq p$. Blue highlights the cases where the tail behavior is constant, so that there is danger of non-robustness. Red highlights the cases where the posterior goes to infinity in the tail, necessarily leading to non-robustness, as this will be shown to be the unique mode.*

LEMMA 3.5. *Assume Assumption 3.1 and Assumption 3.2 hold for each $\mathbf{R}_l$, $1 \leq l \leq p$. If (i) $\gamma_{l_1} \to \infty$ for $1 \leq l_1 \leq p_1$ with $p_1 < p$, (ii) $\gamma_{l_2} \to 0$ for $p_1+1 \leq l_2 \leq p_2$, and (iii) $\gamma_{l_3}$ is bounded between $0$ and $\infty$ for $p_2+1 \leq l_3 \leq p$, then a bound on the tail rate of the marginal posterior of $\boldsymbol{\gamma}$ is*

$$p(\boldsymbol{\gamma} \mid \mathbf{y}^{\mathscr{D}}) \leq C_4 \prod_{l_1=1}^{p_1} |\nu'_{l_1}(\gamma_{l_1})| \left[ \prod_{l_2=p_1+1}^{p_2} \operatorname{tr}\left( \frac{\partial \mathbf{R}}{\partial \gamma_{l_2}} \right)^2 \right]^{1/2},$$

*where $C_4 > 0$ is a positive constant.*

The following corollary is a direct consequence of the above lemma and states that, when the power exponential correlation is used, the only possible infinite mode of the marginal posterior of $\tilde{\boldsymbol{\beta}}$ is at $\tilde{\beta}_l \to 0$ for all $1 \leq l \leq p$.

COROLLARY 3.3. *For the power exponential correlation function, if there is one $l$, $1 \leq l \leq p$, for which $\tilde{\beta}_l > K$ where $K$ is a positive constant, then the marginal posterior of $\tilde{\boldsymbol{\beta}}$ using the standard reference prior (2.7) with $a = 1$ satisfies*

$$p(\tilde{\boldsymbol{\beta}} \mid \mathbf{y}^{\mathscr{D}}) \leq O(1).$$

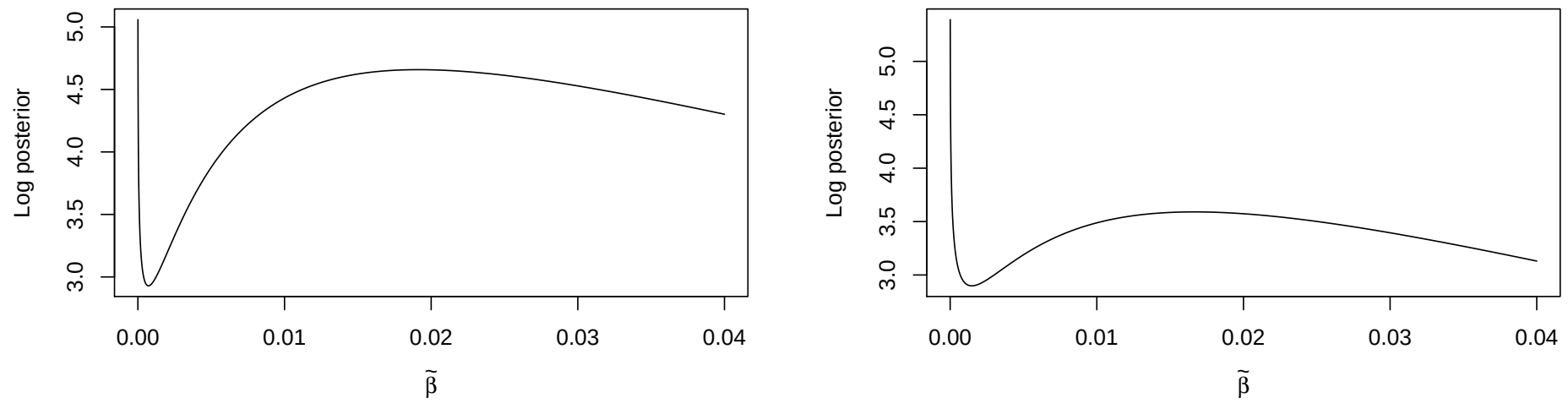

FIG 3. *Examples of the marginal posterior of* $\tilde{\boldsymbol{\beta}}$ *in the power exponential family with* $\alpha = 1.9$*, when emulating the modified Branin function ([9]), which has* $p = 2$ *inputs. Two data sets of size* $n = 20$ *were generated using uniform designs at* $[0,1]^2$ *with* $\mathbf{1} \in \mathcal{C}(\mathbf{h}(\mathbf{x}^{\mathscr{D}}))$*. The black curves are the log marginal posterior of* $\tilde{\beta}$ *arising from setting* $\tilde{\beta}_1 = \tilde{\beta}_2 = \tilde{\beta}$*, and both exhibit infinite posterior density at the mode of 0.*

3.4. *Posterior propriety.* Propriety of the posterior distribution for $\boldsymbol{\gamma}$ (and, hence, for all other parameterizations) is established in the following theorem for general designs, generalizing the theorems in [5] under the isotropic assumption and in [28] for separable designs. For simplicity, we assume $\alpha_1 = \alpha_2 = ... = \alpha_p = \alpha$.

THEOREM 3.2. *When* $\alpha_1 = \alpha_2 = ... = \alpha_p = \alpha$*, the reference prior in (2.7) with* $a = 1$ *results in a proper posterior for GaSP models with the power exponential, spherical, rational quadratic and Matérn correlation functions, under general* $p$*-dimensional designs.*

**4. Robust inference when noise is added to the GaSP model.** Some inputs have little effect on the output of the computer model. Such inputs are called inert inputs ([22]) and are usually not used in building the emulator ([37, 12]). However, when inert inputs are omitted in the emulator, the emulator can no longer be an interpolator at the design points so that the GaSP model is then inappropriate. The common solution is to add a small noise term (sometimes called a nugget) to account for the error, such as $\tilde{y}(\cdot) = y(\cdot) + \epsilon$, where $y(\cdot)$ is the noise-free GaSP and $\epsilon$ is i.i.d. mean-zero Gaussian white noise. This section handles the case where the noise is present in the model. The proofs of the lemmas and theorems in this section are provided in the supplementary materials ([14]).

4.1. *Parameter estimation.* After adding the noise, the covariance function for the new process $\tilde{y}(\cdot)$ can be expressed as

$$\sigma^2 \tilde{c}(\mathbf{x}_l, \mathbf{x}_m) := \sigma^2 \{c(\mathbf{x}_l, \mathbf{x}_m) + \eta \delta_{lm}\}, \tag{4.1}$$

where $\eta$ is defined to be the nugget-variance ratio and $\delta_{lm}$ is a Dirac delta function when $l = m$, i.e., $\delta_{mm} = 1$ and $\delta_{lm} = 0$ if $l \neq m$. Using this parameterization enables marginalization of the likelihood over $\sigma^2$ (c.f., [32]). After adding the noise, the covariance matrix becomes

$$\sigma^2 \tilde{\mathbf{R}} = \sigma^2 (\mathbf{R} + \eta \mathbf{I}_n). \tag{4.2}$$

The reference prior for a real-valued output and isotropic GaSP model with a nugget has been discussed in [32, 19]. Extending it to the GaSP model with multiple range parameters results in the following form:

$$\pi^{\tilde{R}}(\boldsymbol{\theta}, \sigma^2, \boldsymbol{\gamma}, \eta) = \pi^{\tilde{R}}(\boldsymbol{\theta}, \sigma^2)\, \pi^{\tilde{R}}(\boldsymbol{\gamma}, \eta \mid \boldsymbol{\theta}, \sigma^2) \propto \frac{\pi^{\tilde{R}}(\boldsymbol{\gamma}, \eta)}{(\sigma^2)^a}, \tag{4.3}$$

with $\pi^{\tilde{R}}(\boldsymbol{\gamma}, \eta) \propto |\tilde{\mathbf{I}}^*(\boldsymbol{\gamma}, \eta)|^{1/2}$, $\tilde{\mathbf{I}}^*(\cdot)$ the expected Fisher information matrix,

$$\tilde{\mathbf{I}}^*(\boldsymbol{\gamma}, \eta) = \begin{pmatrix} n-q & \mathrm{tr}(\tilde{\boldsymbol{W}}_1) & \mathrm{tr}(\tilde{\mathbf{W}}_2) & \ldots & \mathrm{tr}(\tilde{\mathbf{W}}_{p+1}) \\ & \mathrm{tr}(\tilde{\mathbf{W}}_1^2) & \mathrm{tr}(\tilde{\mathbf{W}}_1 \tilde{\mathbf{W}}_2) & \ldots & \mathrm{tr}(\tilde{\mathbf{W}}_1 \tilde{\mathbf{W}}_{p+1}) \\ & & \mathrm{tr}(\tilde{\mathbf{W}}_2^2) & \ldots & \mathrm{tr}(\tilde{\mathbf{W}}_2 \tilde{\mathbf{W}}_{p+1}) \\ & & & \ddots & \vdots \\ & & & & \mathrm{tr}(\tilde{\mathbf{W}}_{p+1}^2) \end{pmatrix}, \tag{4.4}$$

where $\tilde{\mathbf{W}}_l = \dot{\tilde{\mathbf{R}}}_l \tilde{\mathbf{Q}}$, for $1 \leq l \leq p$, $p$ is the number of range parameters in the correlation matrix $\tilde{\mathbf{R}}$, $\dot{\tilde{\mathbf{R}}}_l$ is the partial derivative of $\tilde{\mathbf{R}}$ with respect to the $l^{th}$ range parameter, and $\tilde{\mathbf{Q}} = \tilde{\mathbf{R}}^{-1} \mathbf{P}_{\tilde{\mathbf{R}}}$ with $\mathbf{P}_{\tilde{\mathbf{R}}} = \mathbf{I}_n - \mathbf{h}(\mathbf{x}^{\mathscr{D}}) \{\mathbf{h}(\mathbf{x}^{\mathscr{D}}) \tilde{\mathbf{R}}^{-1} \mathbf{h}(\mathbf{x}^{\mathscr{D}})\}^{-1} \mathbf{h}(\mathbf{x}^{\mathscr{D}}) \tilde{\mathbf{R}}^{-1}$.

As in the previous sections, one can estimate the nugget and range parameters by their marginal maximum posterior mode,

$$(\hat{\gamma}_1, \ldots \hat{\gamma}_p, \hat{\eta}) = \underset{\gamma_1, \ldots, \gamma_p, \eta}{argmax} \left\{ L(\gamma_1, \ldots, \gamma_p, \eta \mid \mathbf{y}^{\mathscr{D}})\, \pi^{\tilde{\mathbf{R}}}(\gamma_1, \ldots, \gamma_p, \eta) \right\}. \tag{4.5}$$

4.2. *Robustness of the posterior mode.* Note that

$$\tilde{\mathbf{R}} = \mathbf{R}_1 \circ \mathbf{R}_2 \circ \ldots \circ \mathbf{R}_p \circ \mathbf{R}_{p+1},$$

where $\mathbf{R}_{p+1} = \mathbf{1}_n \mathbf{1}_n^T + \eta \mathbf{I}_n$. Also, $\mathbf{R}_{p+1}$ satisfies Assumption 3.2 with $\nu_{p+1}(\eta) = \eta$ and $\omega_{p+1}(\eta) = 0$. Using these facts and in parallel to Lemma 3.3 and Lemma 3.4, the tail rates of the likelihood and the prior for the GaSP with a nugget are given in the following lemmas.

LEMMA 4.1. *If Assumption 3.1 and Assumption 3.2 hold for each of the $\mathbf{R}_l$, $1 \leq l \leq p$ , the marginal likelihood and profile likelihood have the following tail rates.*

*(i) If $\exists l$, $1 \leq l \leq p$, such that $\gamma_l \to 0$, the marginal likelihood and profile likelihood both exist and are greater than zero.*
*(ii) If $\gamma_l \to \infty$ for all $l$, $1 \leq l \leq p$,*

$$L(\boldsymbol{\gamma}, \eta \mid \mathbf{y}^{\mathscr{D}}) = \begin{cases} O\left( (\sum_{l=1}^{p} \nu_l(\gamma_l) + \eta)^{a-1/2} \right), & \mathbf{1}_n \notin \mathcal{C}(\mathbf{h}(\mathbf{x}^{\mathscr{D}})), \\ O\left( (\sum_{l=1}^{p} \nu_l(\gamma_l) + \eta)^{a-1} \right), & \mathbf{1}_n \in \mathcal{C}(\mathbf{h}(\mathbf{x}^{\mathscr{D}})), \end{cases}$$

*and the profile likelihood, in this case, satisfies*

$$L(\boldsymbol{\gamma} \mid \mathbf{y}^{\mathscr{D}}, \hat{\sigma}^2_{MLE}, \hat{\boldsymbol{\theta}}_{MLE}) = O(\sum_{l=1}^{p} \nu_l(\gamma_l) + \eta)^{1/2}.$$

LEMMA 4.2. *If Assumption 3.1 and Assumption 3.2 hold for each of the $\mathbf{R}_l$, $1 \leq l \leq p$, then $\pi^{\tilde{R}}(\gamma, \eta)$ has the following two limiting properties. Here $\boldsymbol{\gamma}_E$ denotes the vector of $\gamma_l$ for all $l \in E$, $E \in \{1, 2, ..., p\}$, and $\boldsymbol{\gamma}_{-E}$ denotes the complementary vector.*

*(i) When $\boldsymbol{\gamma}_E \to \mathbf{0}$ for all $l \in E$, $E \subset \{1, 2, ..., p\}$, then*

$$\pi^R(\boldsymbol{\gamma}) \leq \tilde{C}(\boldsymbol{\gamma}_{-E}) \left[ \prod_{l \in E} \mathrm{tr} \left( \frac{\partial \tilde{\mathbf{R}}}{\partial \gamma_l} \right)^2 \right]^{1/2},$$

*where $\tilde{C}(\boldsymbol{\gamma}_{-E})$ is a constant in $\boldsymbol{\gamma}_E$.*
*(ii) As $\gamma_l \to \infty$ for all $1 \leq l \leq p$ and $\eta \to 0$, if $\mathbf{1} \notin \mathcal{C}(\mathbf{h}(\mathbf{x}^{\mathscr{D}}))$, then*

$$\pi^R(\boldsymbol{\gamma}) \leq \tilde{C}_1 \left| \frac{\prod_{l=1}^{p} \nu_l'(\gamma_l)}{(\sum_{l=1}^{p} \nu_l(\gamma_l) + \eta)^{p+1}} \right|;$$

*further, if $\mathbf{1} \in \mathcal{C}(\mathbf{h}(\mathbf{x}^{\mathscr{D}}))$ and $p \geq 2$,*

$$\pi^R(\boldsymbol{\gamma}) \leq \tilde{C}_2 \left| \frac{\prod_{l=1}^{p} \nu_l'(\gamma_l)}{(\sum_{l=1}^{p} \nu_l(\gamma_l) + \eta)^{p+1}} \right| \left| \sum_{l=1}^{p} \frac{\nu_l^2(\gamma_l) \omega_l'(\gamma_l)}{\nu_l'(\gamma_l) \nu_m(\gamma_m)} \right|,$$

| | $\mathbf{1}_n \in \mathcal{C}(\mathbf{h}(\mathbf{x}^{\mathscr{D}}))$ | | $\mathbf{1}_n \notin \mathcal{C}(\mathbf{h}(\mathbf{x}^{\mathscr{D}}))$ | |
|---|---|---|---|---|
| | $l \in E$, $\gamma_l \to 0$ | $\gamma_l \to \infty$ for all $l$ and $\eta \to 0$ | $l \in E$, $\gamma_l \to 0$ | $\gamma_l \to \infty$ for all $l$ and $\eta \to 0$ |
| Profile Lik | $O(1)$ | $O((\gamma_{(1)}^{-\alpha} + \eta)^{\frac{1}{2}})$ | $O(1)$ | $O((\gamma_{(1)}^{-\alpha} + \eta)^{\frac{1}{2}})$ |
| Marginal Lik | $O(1)$ | $O(1)$ | $O(1)$ | $O((\gamma_{(1)}^{-\alpha} + \eta)^{\frac{1}{2}})$ |
| Post $\boldsymbol{\gamma}$, $p=1$ | $O(\frac{\exp(-C/\gamma^{\alpha})}{\gamma^{(\alpha+1)}})$ | $O(\frac{\gamma^{-2\alpha-1}}{(\gamma^{-\alpha}+\eta)^2})$ | $O(\frac{\exp(-C/\gamma^{\alpha})}{\gamma^{(\alpha+1)}})$ | $O(\frac{\gamma^{-\alpha-1}}{(\gamma^{-\alpha}+\eta)^{3/2}})$ |
| $p \geq 2$ | $O(\prod_{l\in E} \frac{\exp(-C_l/\gamma_l^{\alpha})}{\gamma_l^{(\alpha+1)}})$ | $O(\frac{\prod_{l=1}^{p} \gamma_l^{-\alpha-1} \gamma_{(1)}^{-\alpha}}{(\gamma_{(1)}^{-\alpha}+\eta)^{p+1}})$ | $O(\prod_{l\in E} \frac{\exp(-C_l/\gamma_l^{\alpha})}{\gamma_l^{(\alpha+1)}})$ | $O(\frac{\prod_{l=1}^{p} \gamma_l^{-\alpha-1}}{(\gamma_{(1)}^{-\alpha}+\eta)^{p+1/2}})$ |
| Post $\tilde{\boldsymbol{\beta}}$, $p=1$ | $O(\exp(-\tilde{\beta} C))$ | $O(\frac{\tilde{\beta}}{(\tilde{\beta}+\eta)^2})$ | $O(\exp(-\tilde{\beta} C))$ | $O((\tilde{\beta}+\eta)^{-3/2})$ |
| $p \geq 2$ | $O(\prod_{l\in E} \exp(-\tilde{\beta}_l C_l))$ | $O(\frac{\beta_{(p)}}{(\tilde{\beta}_{(p)}+\eta)^{p+1}})$ | $O(\prod_{l=1}^{p} \exp(-\tilde{\beta}_l) C_l)$ | $O((\tilde{\beta}_{(p)}+\eta)^{-p-1/2})$ |
| Post $\tilde{\boldsymbol{\xi}}$, $p=1$ | $O(\exp(-\exp(\tilde{\xi}) C + \tilde{\xi}))$ | $O(\frac{\exp(2\tilde{\xi})}{(\exp(\tilde{\xi})+\eta)^2})$ | $O(\exp(-\exp(\tilde{\xi}) C))$ | $O(\frac{\exp(\tilde{\xi})}{(\exp(\tilde{\xi})+\eta)^{3/2}})$ |
| $p \geq 2$ | $O(\prod_{l\in E} \exp(-\exp(\tilde{\xi}_l) C_l + \tilde{\xi}_l))$ | $O(\frac{\exp(\sum_{l=1}^{p} \tilde{\xi}_l)\exp(\tilde{\xi}_{(p)})}{(\exp(\tilde{\xi}_{(p)})+\eta)^{p+1}})$ | $O(\prod_{l\in E} \exp(-\exp(\tilde{\xi}_l) C_l))$ | $O(\frac{\exp(\sum_{l=1}^{p} \tilde{\xi}_l)}{(\exp(\tilde{\xi}_{(p)})+\eta)^{p+1/2}})$ |

TABLE 3

*Tail behaviors of the profile likelihood, the marginal likelihood and the posterior distributions for different parameterizations of the power exponential correlation function, using the reference prior in (4.3) with $a = 1$. In the 2nd and 4th columns, $E$ is a nonempty set such that for $l \in E$, $\gamma_l \to 0$ (equivalent to $\tilde{\beta}_l \to \infty$ or $\tilde{\xi}_l \to \infty$), and $C$ and $C_l$ are positive numbers not depending on $\gamma_l \in E$ (or $\tilde{\beta}_l \in E$ or $\tilde{\xi}_l \in E$). In the 3rd and 5th columns, $\gamma_l \to \infty$ (equivalent to $\tilde{\beta}_l \to 0$ or $\tilde{\xi}_l \to -\infty$), for all $1 \leq l \leq p$; in the stated tail rates, $\gamma_{(1)}$ is defined as minimum of the $\gamma_l$, $\tilde{\beta}_{(p)}$ is the largest $\tilde{\beta}_l$, and $\tilde{\xi}_{(p)}$ is the largest $\tilde{\xi}_l$, where $1 \leq l \leq p$. Blue highlights the cases where the tail behavior is constant; red highlights the cases where the posterior goes to infinity in the tail; and green highlights situations in which the rate might go to zero, a constant or infinity, depending on the speed of $\eta$ and $\gamma_l$ to their limits and the choice of the roughness parameter $\alpha$.*

*for every index $m$ between 1 to $p$; if $\mathbf{1} \in \mathcal{C}(\mathbf{h}(\mathbf{x}^{\mathscr{D}}))$ and $p = 1$,*

$$\pi^R(\boldsymbol{\gamma}) \leq \tilde{C}_3 \frac{\nu_1(\gamma_1)|\omega_1'(\gamma_1)|}{(\nu_1(\gamma_1) + \eta)^2},$$

*where $\tilde{C}_1, \tilde{C}_2$ and $\tilde{C}_3$ are positive constants.*

Directly applying Lemma 4.1 and Lemma 4.2 yields the bounds on the tail rates of the marginal posterior under the various parameterizations (and also for the profile and marginal likelihood) in Table 3. For simplicity, we assume $\alpha_1 = \alpha_2 = ... = \alpha_p = \alpha$.

Comparing Table 2 with Table 3, it is clear that addition of the nugget can cause a loss of robustness of the posterior mode for the $(\gamma_1, \gamma_2, ..., \gamma_p, \eta)^T$ and $(\xi_1, \xi_2, ..., \xi_p, \eta)^T$ parameterizations, in certain cases. Luckily, a simple reparameterization of $\eta$, to $\tau = log(\eta)$, with estimation by the corresponding posterior mode, achieves robustness, as shown in the following theorem.

THEOREM 4.1. *When $a = 1$, marginal posterior mode estimation of $(\gamma_1, ..., \gamma_p, \tau)^T$, and $(\xi_1, ..., \xi_p, \tau)^T$, where $\tau = log(\eta)$, is robust for the prod-* *uct form of the power exponential family, spherical, and Matérn correlation functions listed in Table 1, and for the rational quadratic correlation function when* $\alpha > 1/2$*. In addition, marginal posterior mode estimation of* $(\xi_1, ..., \xi_p, \tau)^T$*, for* $1 \leq l \leq p$*, is robust for the rational quadratic correlation function for all* $\alpha > 0$*,* $1 \leq l \leq p$*.*

PROOF. Theorem 4.1 can be proved by verifying Corollary 3.1 and Corollary 3.2, using the results from Lemma 4.1 and Lemma 4.2. □

4.3. *Posterior propriety for the GaSP model with noise.* Propriety of the posterior distribution for $\boldsymbol{\gamma}$ and $\eta$ (and, hence, for all other parameterizations) is established in the following theorem, generalizing the theorems in [32, 19] under the isotropic assumption with a nugget. It can be proved in the same way as Theorem 3.2, so we omit the details. For simplicity, we assume $\alpha_1 = \alpha_2 = ... = \alpha_p = \alpha$.

THEOREM 4.2. *When* $\alpha_1 = \alpha_2 = \ldots = \alpha_p = \alpha$*, the reference prior in (4.3) with* $a = 1$ *results in a proper posterior for the GaSP models with noise, under the power exponential, spherical, rational quadratic and Matérn correlation functions, for general* $p$*-dimensional designs.*

## 5. Numerical results.

5.1. *Comparison criteria.* In this section, we numerically compare the performance of several of the methods discussed above, including the MLE and marginal posterior mode estimation with parameterizations $\boldsymbol{\gamma}$ and $\boldsymbol{\xi}$ (the log inverse of $\boldsymbol{\gamma}$). We do not include the MMLE method or results for the $\tilde{\boldsymbol{\beta}}$ parameterization because of the robustness problems these methods have, as indicated in Table 2 and Table 3. A constant GaSP mean is assumed for all cases, i.e. $h(\mathbf{x}) = 1$, and we use the Matérn correlation with $\alpha = 5/2$ in (2.5) for all methods. Also included are the results produced by the DiceKriging package ([34]), where the Matérn correlation is also the default setting.

We mainly compare the out of sample prediction evaluated by Mean Square Error (MSE). In each simulation, we use $n$ runs, where $n$ is small (typically chosen to be $n \approx 10p$), to build the GaSP emulator, and then record the out-of-sample MSE of $n^* = 10,000$ held-out outputs. This is repeated for $N = 500$ random designs, with the resulting average MSE being reported. The criteria are thus

$$\text{MSE}_j = \frac{1}{n^*}\sum_{i=1}^{n^*}\left(y(\mathbf{x}^*_{ij}) - \hat{y}(\mathbf{x}^*_{ij})\right)^2, \text{ and AvgMSE} = \sum_{\text{j=1}}^{\text{N}} \text{MSE}_\text{j}/\text{N},$$

| | Robust GaSP $\boldsymbol{\xi}$ | Robust GaSP $\boldsymbol{\gamma}$ | MLE | DiceKriging |
|---|---|---|---|---|
| 1-dim Higdon | $1.1 \times 10^{-3}$ | $1.1 \times 10^{-3}$ | $1.2 \times 10^{-3}$ | $1.2 \times 10^{-3}$ |
| 2-dim Branin | $4.7 \times 10^{-7}$ | $4.2 \times 10^{-7}$ | $2.4 \times 10^{-4}$ | $1.8 \times 10^{-3}$ |
| 3-dim D&P | $8.0 \times 10^{-2}$ | $1.5 \times 10^{-1}$ | $8.0 \times 10^{-1}$ | $5.7 \times 10^{-1}$ |
| 4-dim G&L | $4.2 \times 10^{-3}$ | $1.3 \times 10^{-2}$ | $2.8 \times 10^{-2}$ | $4.9 \times 10^{-2}$ |
| 10-dim Linkletter | $1.7 \times 10^{-12}$ | $2.4 \times 10^{-12}$ | $4.8 \times 10^{-5}$ | $5.7 \times 10^{-4}$ |

TABLE 4

*Average MSE of the four estimation procedures for the five experimental functions. From the upper to the lower rows, the sample size is $n = 15, 20, 30, 40$ and 40 for these five functions, respectively. Designs are generated by maxmin LHD. The baseline MSE is 0.52, 36, 52, 0.52, and 0.0044 for these five functions if only the mean of the training output is used for the predictions.*

where $\mathbf{x}^*_{ij}$ is the $i^{th}$ held-out input in the $j^{th}$ design and $\hat{y}(\mathbf{x}^*_{ij})$ is its prediction. To provide a better visual comparison between the methods, we also study the out-of-sample Normalized-RMSE

$$\text{Normalized-RMSE}_j = \sqrt{\sum_{i=1}^{n^*} (y(\mathbf{x}^*_{ij}) - \hat{y}(\mathbf{x}^*_{ij}))^2 / \sum_{i=1}^{n^*} (y(\mathbf{x}^*_{ij}) - \bar{\mathbf{y}}_j)^2},$$

where $\bar{\mathbf{y}}_j$ is the mean of the observed output for the $j^{th}$ experiment, $j = 1, ..., N$. For an effective method, this should range from 0 to 1.

5.2. *GaSP model without a nugget.* We test the following five functions (implemented in [39]):

i. 1 dimensional Higdon function from [17],
   $Y = \sin(2\pi X/10) + 0.2\sin(2\pi X/2.5)$, where $X \in [0, 10]$.
ii. 2 dimensional Branin function from [9],
   $Y = [X_2 - 5.1X_1^2/(4\pi^2) + 5X_1/\pi - 6]^2 + 10[1 - 1/(8\pi)]\cos(X_1) + 10$ where $X_i \in [0, 1]$, for $i = 1, 2$.
iii. 3 dimensional Dette & Pepelyshev function from [7],
   $Y = 4(X_1 - 2 + 8X_2 - 8X_2^2)^2 + (3 - 4X_2)^2 + 16\sqrt{X_3 + 1}(2X_3 - 1)^2$, where $X_i \in [0, 1]$, for $i = 1, 2, 3$.
iv. 4 dimensional modified Gramacy & Lee function from [11],
   $Y = 2\exp\{\sin[0.9^8(X_1 + 0.48)^8]\} + X_2X_3 + X_4$, where $X_i \in [0, 1)$, for $i = 1, 2, 3, 4$.
v. 10 dimensional Linkletter decreasing coefficient function from [22],
   $Y = 0.2(X_1 + X_2/2 + X_3/4 + X_4/8 + X_5/16 + X_6/32 + X_7/64 + X_8/128)$, where $X_i \in [0, 1]$, for $i = 1$ to 10. Only the first eight inputs are effective.

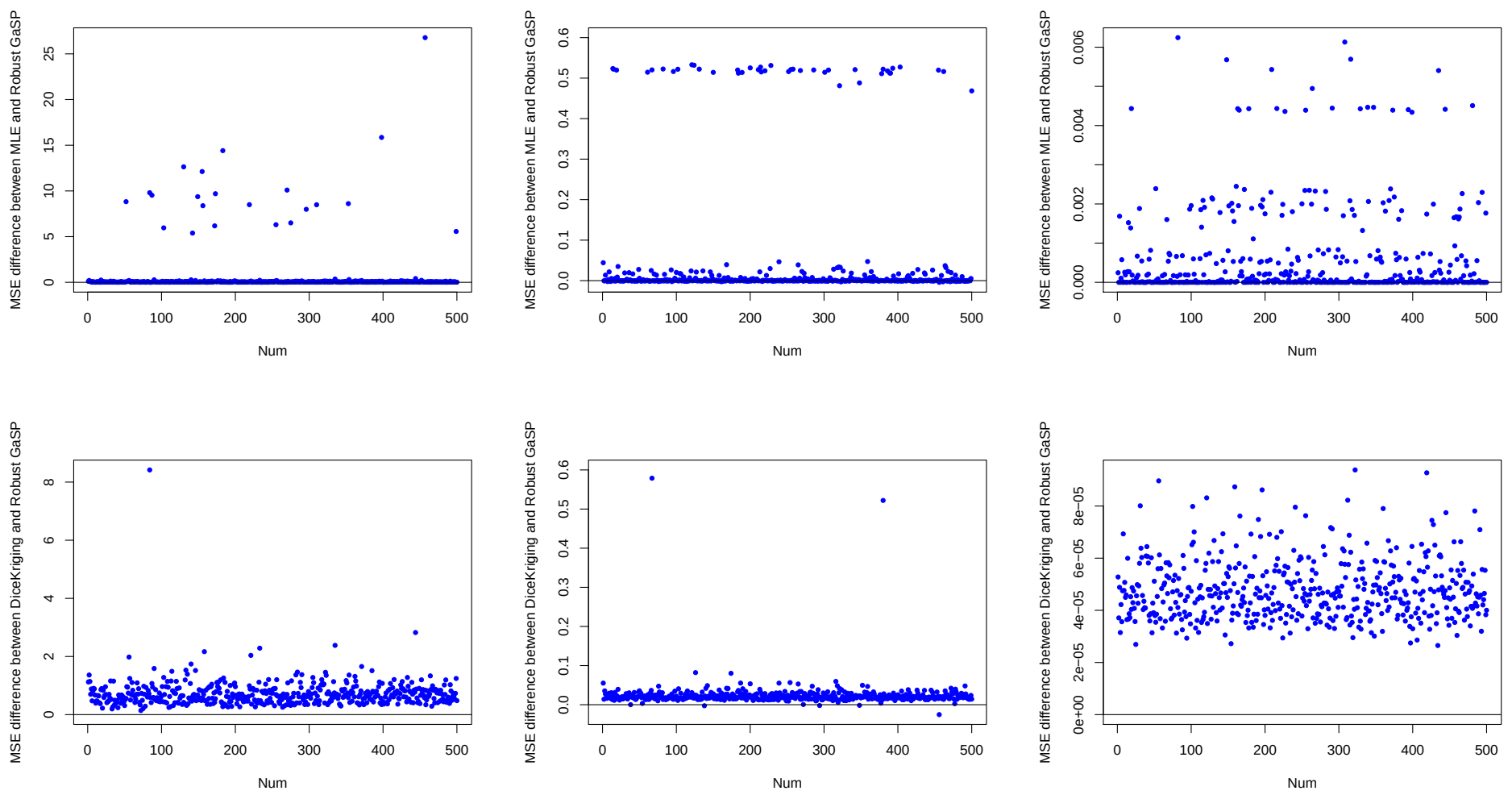

FIG 4. *Plots of $MSE$ difference for each of $N = 500$ designs for the Dette and Pepelyshev function (left panels), Gramacy and Lee function (middle panels), and Linkletter decreasing coefficient function (right panels). The $MSE$ for the MLE GaSP minus the $MSE$ for the robust GaSP under the $\boldsymbol{\xi}$ parameterization is plotted in the first row, and the $MSE$ for DiceKriging minus the $MSE$ for the robust GaSP under the $\boldsymbol{\xi}$ parameterization is plotted in the second row.*

The average MSEs of the four estimation methods for the five functions are shown in Table 4. The robust GaSP methods were implemented using [13] and they clearly outperformed the MLE and DiceKriging, with the $\boldsymbol{\xi}$ parameterization yielding the best performance for most of the cases. Note that all methods used the same GaSP prediction equations; the only difference was in the estimates of the correlation parameters.

The first row in Figure 4 gives the difference of $\text{MSE}_j$ of prediction, for each of 500 designs $j$ (for functions iii, iv and v), between the MLE GaSP and the robust GaSP under the $\boldsymbol{\xi}$ parameterization. Note that, for a significant proportion of the designs, the MLE GaSP is much worse than the robust GaSP. In these cases, the MLE GaSP estimate yields a covariance matrix that is close to $\hat{\mathbf{R}} \approx \mathbf{I}_n$, so that the prediction degenerated to the fitted mean with impulse functions at the observed values of the inputs.

The second row in Figure 4 gives the difference of $\text{MSE}_j$ of prediction, for each of 500 designs $j$, between the DiceKriging GaSP and the robust GaSP under the $\boldsymbol{\xi}$ parameterization. The DiceKriging package uses a number of techniques to avoid unstable prediction of the correlation parameters ([34]) and is more stable than the MLE (without any constraints), as can be seen by a comparison of the upper panels and the lower panels of Figure 4 (the $y$-axis scales are considerably smaller for DiceKriging). Clearly, however, DiceKriging produces inferior correlation parameter estimates than does the robust GaSP in virtually all of the design cases for the three functions in Figure 4; indeed, only for few design choices for the Gramacy & Lee function does DiceKriging produce better predictions than the robust GaSP.

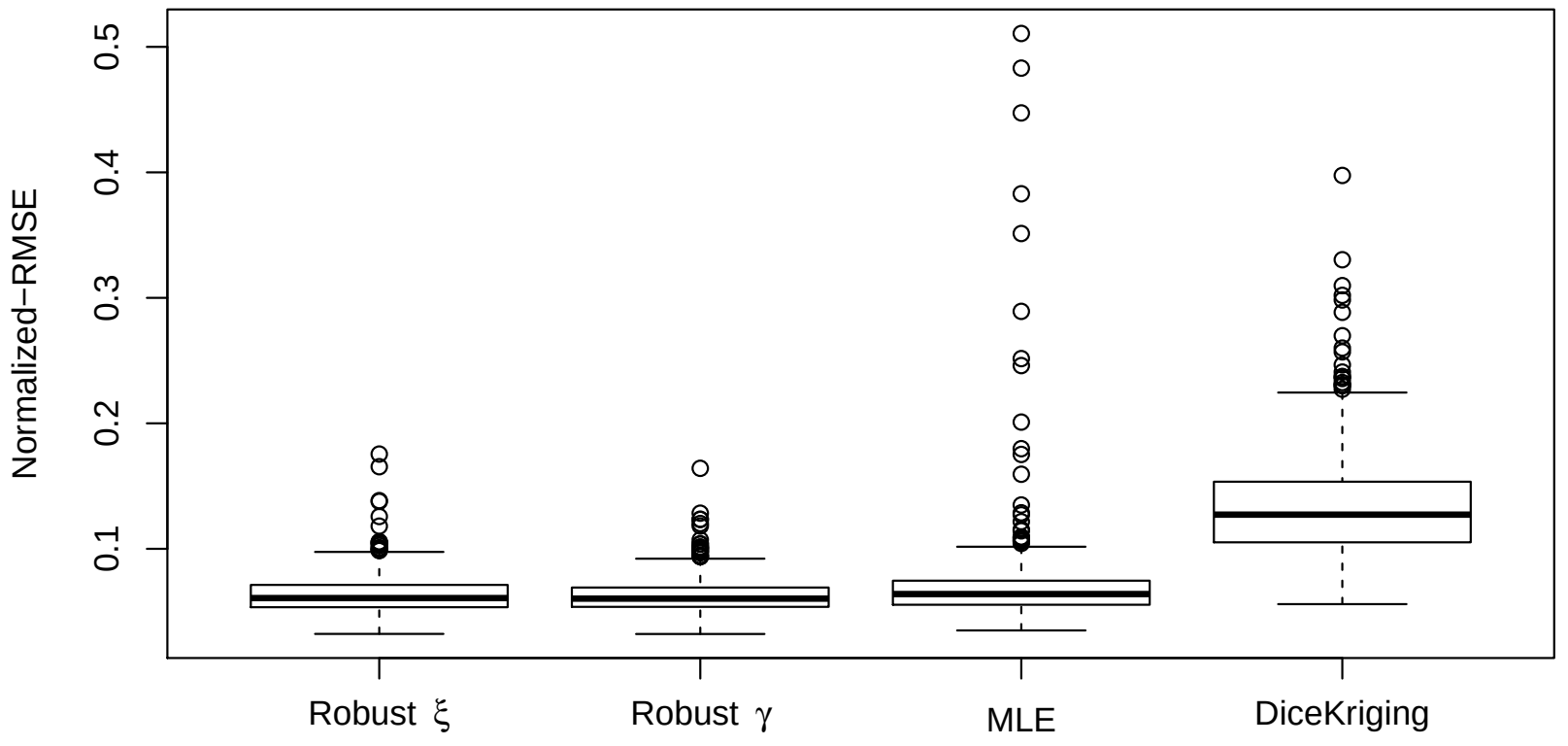

FIG 5. *Boxplots, for the four estimation methods, of the Normalized RMSE for prediction of the Borehole function, based on* $n = 25$ *design points to build the emulator and averaging over* $N = 500$ *different designs generated from a Maximin LHD design. The average baseline MSE is 2079.36, using only the mean for prediction. The Average MSE for the 4 methods (from the left to the right) are 9.07, 8,72, 14.77 and 41.29.*

5.3. *GaSP model with noise.* The borehole function models water flow through a borehole ([1, 24]) and is given by

$$Y = \frac{2\pi T_u(H_u - H_l)}{\ln(r/r_\omega)\{1 + 2LT_u/[\ln(r/r_\omega)r_\omega^2 K_\omega] + T_u/T_l\}},$$

where $r_\omega, r, T_u, H_u, T_l, H_l, L, K_\omega$ are the 8 inputs. The inputs $r$, $T_u$ and $T_l$ barely affect the output (as clearly shown in Figure S1 of the supplementary materials [14], where we draw plots of the borehole function by fixing seven of the inputs and varying one), and this holds globally over the input space. We thus only use the remaining five influential inputs to build the GaSP model, and then add a nugget to account for the error.

The results of Normalized-RMSE for the borehole function are shown in Figure 5. The average MSE of the GaSP with parameters estimated by MLE is 14.77, which is worse than the robust GaSP with the $\boldsymbol{\xi}$ or $\boldsymbol{\gamma}$ parameterization, whose average MSEs are 9.07 and 8.72, respectively. This is because $\tilde{\mathbf{R}} \to (1+\hat{\eta})\mathbf{I}_n$ for the MLE in several of the cases. Although the nugget might stabilize the computation when $\mathbf{R} \approx \mathbf{1}_n\mathbf{1}_n^T$, it cannot help when $\mathbf{R}$ becomes nearly proportional to $\mathbf{I}_n$. In contrast, the robust GaSP, with a good parameterization, prevents these bad cases from materializing.

DiceKriging is worse than the robust GaSP in almost all cases, in terms of Normalized-RMSE. The average Normalized-RMSE for robust GaSP, with the $\boldsymbol{\xi}$ or $\boldsymbol{\gamma}$ parameterization, are around 0.064 and 0.063 respectively, both of which are quite small, considering that only $n = 25$ observations were utilized to build the emulators.

**6. Concluding remarks.** We have introduced the robust GaSP for computer model emulation, namely marginal posterior mode estimation of the emulating Gaussian process correlation parameters, using the reference prior and certain parameterizations. This emulation methodology was shown to be robust for a wide class of correlation functions, whereas a number of alternative methods were shown to be non-robust. Robustness here means that the estimates of the correlation parameters avoid two possibly severe problems that can happen: the estimated correlation matrix could be nearly singular or could nearly equal the identity matrix. We also proved posterior propriety, under the reference prior, for general multi-dimensional designs.

The current study of the tail behavior of the likelihoods and posteriors can be extended to the situation where both the roughness and range parameters are unknown and to a more general class of correlation functions. In addition, the results hold when the inputs are divided into $k < p$ groups, and an isotropic correlation function is assumed for each group. The results about tail behavior, given here, were for a finite number of observations and it would be interesting to understand the tail behavior as the number of observations goes to infinity.

**Appendix: The problem of design singularity for power exponential correlation and $\boldsymbol{\alpha = 2}$.** Consider a single input and an equally-spaced design on [0,1], with $n = 10$, so that the inputs are $x_i^{\mathscr{D}} = (i-1)/(n-1)$, $i = 1, \cdots, n$. Suppose one uses the power exponential correlation function in Table 1 with roughness parameter $\alpha = 2$. Denote the "design correlation" matrix as $\mathbf{R}^0$ with the $(i,j)$ entry $|x_i^{\mathscr{D}} - x_j^{\mathscr{D}}|^2$, $1 \leq i,j \leq n$. The condition number of $\mathbf{R}^0$ is larger than $10^{16}$. $\mathbf{R}$ in this case is also ill-conditioned with a small range parameter $\gamma$, e.g., $\gamma = 1$. Although $\mathbf{R}$ is quite far away from $\mathbf{1}_n\mathbf{1}_n^T$, it is near singular and becomes almost non-invertible when $n \geq 15$.

This type of singularity is reported in the literature (c.f., [29]). When $\mathbf{R}^0$ is ill-conditioned, then usually $\mathbf{R}$ is ill-conditioned even if $\mathbf{R}$ is far away from $\mathbf{1}_n\mathbf{1}_n^T$. Clearly, this type of matrix singularity is related to the choice of roughness parameters $\alpha$, but not related to the range parameters $\gamma$.

One remedy for design singularity is to replace Gaussian covariance by Matérn covariance, or simply choose the range parameter $\alpha < 2$ in power exponential correlation as in [4, 37]. This type of singularity is a separate problem from what we considered, and can be avoided by a pre-experimental check of the design correlation matrix.

**Acknowledgements.** This research was supported by NSF grants DMS-1007773, DMS-1228317, EAR-1331353, DMS-1407775, and EAR 1521855. The research of Mengyang Gu was part of his PhD thesis at Duke University. The authors thank the editor, the associate editor and two referees for their comments that substantially improved the article.

**References.**

[1] An, J. and Owen, A. (2001). Quasi-regression. *Journal of complexity* **17** 588–607.

[2] Andrianakis, I. and Challenor, P. G. (2012). The effect of the nugget on Gaussian process emulators of computer models. *Computational Statistics & Data Analysis* **56** 4215–4228.

[3] Bayarri, M., Berger, J., Cafeo, J., Garcia-Donato, G., Liu, F., Palomo, J., Parthasarathy, R., Paulo, R., Sacks, J. and Walsh, D. (2007). Computer model validation with functional output. *The Annals of Statistics* **35** 1874–1906.

[4] Bayarri, M. J., Berger, J. O., Calder, E. S., Dalbey, K., Lunagomez, S., Patra, A. K., Pitman, E. B., Spillerh, E. T. and Wolperti, R. L. (2009). Using statistical and computer models to quantify volcanic hazards. *Technometrics* **51** 402-413.

[5] Berger, J. O., De Oliveira, V. and Sansó, B. (2001). Objective Bayesian analysis of spatially correlated data. *Journal of the American Statistical Association* **96** 1361–1374.

[6] De Oliveira, V. (2007). Objective Bayesian analysis of spatial data with measurement error. *Canadian Journal of Statistics* **35** 283–301.

[7] Dette, H. and Pepelyshev, A. (2010). Generalized latin hypercube design for computer experiments. *Technometrics* **52**.

[8] Diggle, P. and Ribeiro, P. (2007). *Model-based Geostatistics*. Springer.

[9] Dixon, L. (1978). The global optimization problem: an introduction. *Towards Global Optimiation 2* 1–15.

[10] Gelfand, A. E., Diggle, P., Guttorp, P. and Fuentes, M. (2010). *Handbook of spatial statistics*. CRC Press.

[11] Gramacy, R. B. and Lee, H. K. (2009). Adaptive design and analysis of supercomputer experiments. *Technometrics* **51** 130–145.

[12] Gu, M., Berger, J. O. et al. (2016). Parallel partial Gaussian process emulation for computer models with massive output. *The Annals of Applied Statistics* **10** 1317–1347.

[13] Gu, M., Palomo, J. and Berger, J. (2016). RobustGaSP: Robust Gaussian Stochastic Process Emulation R package version 0.5.

[14] Gu, M., Wang, X. and Berger, J. O. (2017). Supplement to "Robust Gaussian Stochastic Process Emulation".

[15] Handcock, M. S. and Stein, M. L. (1993). A Bayesian analysis of kriging. *Technometrics* **35** 403–410.

[16] Handcock, M. S. and Wallis, J. R. (1994). An approach to statistical spatial-temporal modeling of meteorological fields. *Journal of the American Statistical Association* **89** 368–378.

[17] Higdon, D. et al. (2002). Space and space-time modeling using process convolutions. *Quantitative methods for current environmental issues* **37–56**.

[18] Kazianka, H. (2013). Objective bayesian analysis of geometrically anisotropic spatial data. *Journal of Agricultural, Biological, and Environmental Statistics* **18** 514–537.

[19] Kazianka, H. and Pilz, J. (2012). Objective Bayesian analysis of spatial data with uncertain nugget and range parameters. *Canadian Journal of Statistics* **40** 304–327.

[20] Kennedy, M. C. and O'Hagan, A. (2001). Bayesian calibration of computer models. *Journal of the Royal Statistical Society: Series B (Statistical Methodology)* **63** 425–464.

[21] Li, R. and Sudjianto, A. (2005). Analysis of computer experiments using penalized likelihood in Gaussian Kriging models. *Technometrics* **47**.

[22] Linkletter, C., Bingham, D., Hengartner, N., Higdon, D. and Kenny, Q. Y. (2006). Variable selection for Gaussian process models in computer experiments. *Technometrics* **48** 478-490.

[23] Lopes, D. (2011). Development and implementation of Bayesian computer model emulators PhD thesis, Duke University.

[24] Morris, M. D., Mitchell, T. J. and Ylvisaker, D. (1993). Bayesian design and analysis of computer experiments: use of derivatives in surface prediction. *Technometrics* **35** 243–255.

[25] Oakley, J. (1999). Bayesian uncertainty analysis for complex computer codes PhD thesis, University of Sheffield.

[26] Oakley, J. and O'Hagan, A. (2002). Bayesian inference for the uncertainty distribution of computer model outputs. *Biometrika* **89** 769–784.

[27] Paciorek, C. J. and Schervish, M. J. (2006). Spatial modelling using a new class of nonstationary covariance functions. *Environmetrics* **17** 483–506.

[28] Paulo, R. (2005). Default priors for Gaussian processes. *The Annals of statistics* **33** 556–582.

[29] Peng, C.-Y. and Wu, C. J. (2014). On the choice of nugget in kriging modeling for deterministic computer experiments. *Journal of Computational and Graphical Statistics* **23** 151–168.

[30] Qian, P. Z. G., Wu, H. and Wu, C. F. J. (2008). Gaussian Process Models for Computer Experiments With Qualitative and Quantitative Factors. *Technometrics* **50** 383-396.

[31] Ranjan, H. R. P. and Karsten, R. (2011). A Computationally Stable Approach to Gaussian Process Interpolation of Deterministic Computer Simulation Data. *Technometrics* **53** 366 - 378.

[32] Ren, C., Sun, D. and He, C. (2012). Objective Bayesian analysis for a spatial model with nugget effects. *Journal of Statistical Planning and Inference* **142** 1933–1946.

[33] Ren, C., Sun, D. and Sahu, S. K. (2013). Objective Bayesian analysis of spatial models with separable correlation functions. *Canadian Journal of Statistics* **41** 488–507.

[34] Roustant, O., Ginsbourger, D. and Deville, Y. (2012). DiceKriging, DiceOptim: Two R Packages for the Analysis of Computer Experiments by Kriging-Based

Metamodeling and Optimization. *Journal of Statistical Software* **51** 1–55.

[35] Sacks, J., Welch, W. J., Mitchell, T. J., Wynn, H. P. et al. (1989). Design and analysis of computer experiments. *Statistical science* **4** 409–423.

[36] Santner, T. J., Williams, B. J. and Notz, W. I. (2003). *The design and analysis of computer experiments*. Springer Science & Business Media.

[37] Spiller, E. T., Bayarri, M., Berger, J. O., Calder, E. S., Patra, A. K., Pitman, E. B. and Wolpert, R. L. (2014). Automating emulator construction for geophysical hazard maps. *SIAM/ASA Journal on Uncertainty Quantification* **2** 126–152.

[38] Stein, M. L. (2012). *Interpolation of spatial data: some theory for kriging*. Springer Science & Business Media.

[39] Surjanovic, S. and Bingham, D. Virtual Library of Simulation Experiments: Test Functions and Datasets. Retrieved June 26, 2017, from http://www.sfu.ca/~ssurjano.

[40] Zhang, H. (2004). Inconsistent estimation and asymptotically equal interpolations in model-based geostatistics. *Journal of the American Statistical Association* **99** 250–261.

[41] Zhang, H. and Zimmerman, D. L. (2005). Towards reconciling two asymptotic frameworks in spatial statistics. *Biometrika* **92** 921–936.

[42] Zimmerman, D. L. (1993). Another look at anisotropy in geostatistics. *Mathematical Geology* **25** 453–470.

M. Gu
Department of Applied Mathematics and Statistics
Johns Hopkins University,
3400 North Charles Street, Whitehead Hall 100,
Baltimore, Maryland, 21218-2608
USA
E-mail: mgu6@jhu.edu

X. Wang
Department of Statistics
University of Connecticut,
215 Glenbrook Rd. U-4120
Storrs, CT 06269-4120
USA
E-mail: xiaojing.wang@uconn.edu

J. O. Berger
Department of Statistical Science
Duke University, P.O. Box 90251
Durham, North Carolina 27708-0251
USA
E-mail: berger@stat.duke.edu

# SUPPLEMENTARY MATERIALS: ROBUST GAUSSIAN STOCHASTIC PROCESS EMULATION

All the formulas in this supplementary materials are cross-referenced in the main body of the article.

**S1. Proofs for Section 3.1.** The logarithm of the profile likelihood and marginal likelihood have the following forms, respectively.

(S1)

$$\begin{aligned}&\log L(\rho \mid \mathbf{y}^{\mathscr{D}}, \hat{\sigma}^2_{MLE}, \hat{\boldsymbol{\theta}}_{MLE})\\ &= C_0 + \frac{1}{2}\log(1-\rho^2) - \frac{n}{2}\log\left\{\sum_{i=1}^{n} y_i^2 - 2\rho\left(\sum_{i=1}^{n-1} y_i y_{i+1}\right) + \rho^2\left(\sum_{i=2}^{n-1} y_i^2\right)\right.\\ &\quad\left. - \frac{1-\rho}{n-(n-2)\rho}\left[\sum_{i=1}^{n}\sum_{j=1}^{n} y_i y_j - 2\rho\left(\sum_{i=1}^{n}\sum_{j=2}^{n-1} y_i y_j\right) + \rho^2\left(\sum_{i=2}^{n-1}\sum_{j=2}^{n-1} y_i y_j\right)\right]\right\}.\end{aligned}$$

(S2)

$$\begin{aligned}&\log L(\rho \mid \mathbf{y}^{\mathscr{D}})\\ &= \tilde{C}_0 - \frac{1}{2}\log\left(\frac{n-(n-2)\rho}{1+\rho}\right) - \frac{n-1}{2}\log\left\{\sum_{i=1}^{n} y_i^2 - 2\rho\left(\sum_{i=1}^{n-1} y_i y_{i+1}\right) + \rho^2\left(\sum_{i=2}^{n-1} y_i^2\right)\right.\\ &\quad\left. - \frac{1-\rho}{n-(n-2)\rho}\left[\sum_{i=1}^{n}\sum_{j=1}^{n} y_i y_j - 2\rho\left(\sum_{i=1}^{n}\sum_{j=2}^{n-1} y_i y_j\right) + \rho^2\left(\sum_{i=2}^{n-1}\sum_{j=2}^{n-1} y_i y_j\right)\right]\right\}.\end{aligned}$$

where $C_0$ and $\tilde{C}_0$ are two constants unrelated to $\mathbf{y}^{\mathscr{D}}$ and $\rho$. Equation (S2) is shown in [3], and we show the validity of Equation (S1) as follows. Note that the correlation matrix is

$$\mathbf{R} = \begin{pmatrix} 1 & \rho & \rho^2 & \rho^3 & \dots & \rho^{n-1} \\ \rho & 1 & \rho & \rho^2 & \dots & \rho^{n-2} \\ \rho^2 & \rho & 1 & \rho & \dots & \rho^{n-3} \\ \vdots & \vdots & \vdots & \ddots & \vdots & \vdots \\ \rho^{n-2} & \rho^{n-3} & \rho^{n-4} & \dots & 1 & \rho \\ \rho^{n-1} & \rho^{n-2} & \rho^{n-3} & \dots & \rho & 1 \end{pmatrix},$$

and the inverse correlation matrix is

$$\mathbf{R}^{-1} = \frac{1}{1-\rho^2}\begin{pmatrix} 1 & -\rho & 0 & 0 & \dots & 0 \\ -\rho & 1+\rho^2 & -\rho & 0 & \dots & 0 \\ 0 & -\rho & 1+\rho^2 & -\rho & \dots & 0 \\ \vdots & \vdots & \vdots & \ddots & \vdots & \vdots \\ 0 & 0 & \dots & -\rho & 1+\rho^2 & -\rho \\ 0 & 0 & \dots & 0 & -\rho & 1 \end{pmatrix}.$$

Direct computation yields

$$\mathbf{1}_n^T \mathbf{R}^{-1} \mathbf{1}_n = \frac{n-(n-2)\rho}{1+\rho},$$
$$(\mathbf{y}^{\mathscr{D}})^T \mathbf{R}^{-1} (\mathbf{y}^{\mathscr{D}}) = \frac{\sum_{i=1}^n y_i^2 - 2\rho \sum_{i=1}^{n-1} y_i y_{i+1} + \rho^2 \sum_{i=2}^{n-1} y_i^2}{1-\rho^2},$$
$$|\mathbf{R}| = (1-\rho^2)^{n-1},$$

from which the Equation (S1) follows.

As $\rho \to 0$ and $\rho \to 1$, the limiting values for the log-marginal likelihood and log-profile likelihood are given in Table S1. If the mode for either log likelihood is at (or near) zero, the emulator will degenerate, as in Figure 1. From Table S1, it is clear that the log likelihoods will not typically go to zero as $\rho \to 0$.

| | when $\rho \to 0$ | when $\rho \to 1$ |
|---|---|---|
| $\log L(\rho \mid \mathbf{y}^{\mathscr{D}})$ | $-\frac{n-1}{2}\log\{\sum_{i=1}^n (y_i - \bar{y})^2\}$ | $-\frac{n-1}{2}\log\{\sum_{i=1}^{n-1} (y_{i+1} - y_i)^2\}$ |
| $\log L(\rho \mid \mathbf{y}^{\mathscr{D}}, \hat{\sigma}^2_{MLE}, \hat{\boldsymbol{\theta}}_{MLE})$ | $-\frac{n}{2}\log\{\sum_{i=1}^n (y_i - \bar{y})^2\}$ | $-\infty$ |

Table S1

*The tail behaviors of the log-marginal likelihood and log-profile likelihood (up to the normalizing constants unrelated to $\rho$ and $\mathbf{y}^{\mathscr{D}}$).*

To prove Lemma 3.1, the following quantities are needed:

$$a = \sum_{i=1}^n y_i^2, \quad b = \sum_{i=1}^{n-1} y_i y_{i+1}, \quad c = \sum_{i=2}^{n-1} y_i^2,$$
$$d = \sum_{i=1}^n \sum_{j=1}^n y_i y_j, \quad e = \sum_{i=1}^n \sum_{j=2}^{n-1} y_i y_j, \quad f = \sum_{i=2}^{n-1} \sum_{j=2}^{n-1} y_i y_j,$$

and

$$U = a - 2\rho b + \rho^2 c - \frac{1-\rho}{n-(n-2)\rho}\left[d - 2\rho e + \rho^2 f\right].$$

The following lemma is also needed for the proof of Lemma 3.1.

LEMMA S1.1. *We have*

$$na \geq d, \tag{S3}$$

$$(n-2)c \geq f, \tag{S4}$$

$$c[n-(n-2)\rho] > (1-\rho)f. \tag{S5}$$

PROOF. The first two inequalities are obvious. As

$$c[n-(n-2)\rho]-(1-\rho)f = (n-2)c\frac{[n-(n-2)\rho]}{n-2}-(1-\rho)f,$$

the third inequality follows from $(n-2)c \geq f$ and $\frac{[n-(n-2)\rho]}{n-2}-(1-\rho) = \frac{2}{n-2} > 0$. □

PROOF OF LEMMA 3.1. One only needs to prove that the log-profile likelihood in Equation (S1) decreases with $\rho$ for all $0 \leq \rho \leq 1$. The derivative for Equation (S1) is

$$\frac{\partial \mathrm{log} L(\rho \mid \mathbf{y}^{\mathscr{D}}, \hat{\sigma}^2_{MLE}, \hat{\boldsymbol{\theta}}_{MLE})}{\partial \rho} \propto -\frac{\rho}{1-\rho^2}-\frac{n}{2U}\frac{\partial U}{\partial \rho}, \tag{S6}$$

in which the second term can be written as

$$\frac{n}{2U}\frac{\partial U}{\partial \rho}$$
$$=\frac{n\{(c\rho-b)[n-(n-2)\rho]^2+(d-2e\rho+f\rho^2)-(1-\rho)[n-(n-2)\rho](-e+\rho f)\}}{[n-(n-2)\rho]\{(a-2b\rho+c\rho^2)[n-(n-2)\rho]-(1-\rho)(d-2e\rho+f\rho^2)\}}.$$

(a) To show necessity, note first that, as $\rho \to 1$, $logL(\mathbf{y}^{\mathscr{D}}|\hat{\sigma}^2_{MLE}, \hat{\boldsymbol{\theta}}_{MLE}, \rho) \to -\infty$. Thus a necessary condition is that

$$\lim_{\rho \to 0^+} \frac{\partial \mathrm{log} L(\rho \mid \mathbf{y}^{\mathscr{D}}, \hat{\sigma}^2_{MLE}, \hat{\boldsymbol{\theta}}_{MLE})}{\partial \rho} \leq 0,$$

which implies

$$n^2 b - d - ne \leq 0. \tag{S7}$$

(b) To show sufficiency, we separately discuss $e \leq 0$ and $e > 0$.

(b1) When $e \leq 0$, noticing that $U = S^2/(1-\rho^2) > 0$, a sufficient condition is $\frac{\partial U}{\partial \rho} > 0$ or, equivalently, that the numerator of the second term of Equation (S6) is positive. First considering the terms related to $c$ and $f$ and applying Inequality (S5) yields

$$c\rho[n-(n-2)\rho]^2-(1-\rho)\rho[n-(n-2)\rho]f > 0.$$

For the remaining terms,

$$\begin{aligned} n[n-(n-2)\rho]^2 \Bigg\{ & -b + \frac{d}{[n-(n-2)\rho]^2} \\ & + e\left( \frac{-2\rho}{[n-(n-2)\rho]^2} + \frac{1-\rho}{[n-(n-2)\rho]} \right) \Bigg\} \\ & \geq n[n-(n-2)\rho]^2 \{-b + \frac{d}{n^2} + h(\rho)e\}, \end{aligned}$$

where

$$h(\rho) = \frac{-2\rho}{[n-(n-2)\rho]^2} + \frac{1-\rho}{[n-(n-2)\rho]}.$$

It is easy to show $h(\rho)$ is decreasing monotonically with $h(0) = \frac{1}{n} > 0$ and $h(1) = -\frac{1}{2} < 0$. Since $e \leq 0$, we thus have $h(\rho)e \geq h(0)e = \frac{e}{n}$. Thus a sufficient condition is $-b + \frac{d}{n^2} + h(\rho)e \geq -b + \frac{d}{n^2} + h(0)e = -b + \frac{d}{n^2} + \frac{e}{n} \geq 0$, which is equivalent to

$$b - \frac{d}{n^2} - \frac{e}{n} \leq 0.$$

(b2) We show that $b - \frac{d}{n^2} - \frac{e}{n} \leq 0$ is a sufficient condition for $e > 0$ as follows. First, $b - \frac{d}{n^2} - \frac{e}{n} \leq 0$ is equivalent to

$$-b \geq \frac{-d - ne}{n^2}. \tag{S8}$$

Second,

$$e = \sqrt{df} \leq \frac{\lambda d}{2} + \frac{f}{2\lambda}, \tag{S9}$$

for any $\lambda > 0$. After ignoring the constant $n$, the numerator of the second term of Equation (S6) is

$$
\begin{aligned}
&(c\rho-b)[n-(n-2)\rho]^2+(d-2e\rho+f\rho^2)\\
&\quad-(1-\rho)[n-(n-2)\rho](-e+\rho f)\\
=&c\rho\frac{n-2}{n-2}[n-(n-2)\rho]^2-(1-\rho)[n-(n-2)\rho]\rho f\\
&\quad-b[n-(n-2)\rho]^2+(d-2e\rho+f\rho^2)+(1-\rho)[n-(n-2)\rho]e\\
\geq&\rho f\frac{[n-(n-2)\rho]^2}{n-2}-(1-\rho)[n-(n-2)\rho]\rho f-b[n-(n-2)\rho]^2\\
&\quad+(d-2e\rho+f\rho^2)+(1-\rho)[n-(n-2)\rho]e\\
=&\rho f\frac{2}{n-2}[n-(n-2)\rho]-b[n-(n-2)\rho]^2\\
&\quad+(d-2e\rho+f\rho^2)+(1-\rho)[n-(n-2)\rho]e\\
\geq&\rho f\frac{2}{n-2}[n-(n-2)\rho]+\frac{-d-ne}{n^2}[n-(n-2)\rho]^2\\
&\quad+(d-2e\rho+f\rho^2)+(1-\rho)[n-(n-2)\rho]e\\
=&\left(-\frac{[n-(n-2)\rho]^2}{n^2}+1\right)d-\left(2-\frac{n-2}{n}\rho\right)(2\rho e)\\
&\quad+\left(\rho\frac{2}{n-2}[n-(n-2)\rho]+\rho^2\right)f\\
\geq&\left(-\frac{[n-(n-2)\rho]^2}{n^2}+1\right)d-\left(2-\frac{n-2}{n}\rho\right)\rho\left(\lambda d+\frac{f}{\lambda}\right)\\
&\quad+\left(\rho\frac{2}{n-2}[n-(n-2)\rho]+\rho^2\right)f\\
=&\left(-\frac{[n-(n-2)\rho]^2}{n^2}+1-(2-\frac{n-2}{n}\rho)\rho\lambda\right)d\\
&\quad+\left(\rho\frac{2}{n-2}[n-(n-2)\rho]+\rho^2-(2-\frac{n-2}{n}\rho)\frac{\rho}{\lambda}\right)f\,.
\end{aligned}
$$

The first inequality is from the result in Lemma S1.1. The second inequality follows from Equation (S8). The third inequality is from Equation (S9) and the fact that $2-\frac{n-2}{n}\rho\geq 0$. Finally, putting $\lambda=\frac{n-2}{n}$ into the last equation, the coefficients of $d$ and $f$ are as follows.

- ∗ The coefficient of $d$ is $-\frac{[n-(n-2)\rho]^2}{n^2}+1-(2-\frac{n-2}{n}\rho)\rho\frac{n-2}{n}=0$;
- ∗ The coefficient of $f$ is $\rho\frac{2}{n-2}[n-(n-2)\rho]+\rho^2-(2-\frac{n-2}{n}\rho)\rho\frac{n}{n-2}=0$;

from which the proof is complete.

□

**S2. Proofs for Section 3.3.** Here are some needed facts.

Fact 1. (Schur Product Theorem) If $\mathbf{A}_1$ and $\mathbf{A}_2$ are both positive semidefinite matrices (i.o, $\mathbf{A}_1 \geq 0$ and $\mathbf{A}_2 \geq 0$)

$$\mathbf{A}_1 \circ \mathbf{A}_2 \geq 0.$$

Fact 2. For positive semidefinite matrices $\mathbf{A}_1$, $\mathbf{A}_2$, $\mathbf{A}_3$ and $\mathbf{A}_4$, if $\mathbf{A}_1 \geq \mathbf{A}_2$, $\mathbf{A}_3 \geq \mathbf{A}_4$,

$$\mathbf{A}_1 \circ \mathbf{A}_3 \geq \mathbf{A}_2 \circ \mathbf{A}_4.$$

Fact 3. If $\mathbf{A}$ be a nonsingular matrix and $\boldsymbol{u}$ and $\boldsymbol{v}$ are two vectors, then

$$|\mathbf{A} + \boldsymbol{u}\boldsymbol{v}'| = |\mathbf{A}|(1 + \boldsymbol{v}'\mathbf{A}^{-1}\boldsymbol{u}). \tag{S10}$$

Further, if $\boldsymbol{v}'\mathbf{A}^{-1}\boldsymbol{u} \neq -1$, then $\mathbf{A} + \boldsymbol{u}\boldsymbol{v}'$ is nonsingular and

$$(\mathbf{A} + \boldsymbol{u}\boldsymbol{v}')^{-1} = \mathbf{A}^{-1} - \frac{(\mathbf{A}^{-1}\boldsymbol{u})(\boldsymbol{v}'\mathbf{A}^{-1})}{1 + \boldsymbol{v}'\mathbf{A}^{-1}\boldsymbol{u}}, \tag{S11}$$

which is called the Sherman-Morrison Formula.

Fact 4. $\mathbf{D}_l$ defined in Assumption 3.2 has $n-1$ positive eigenvalues and one negative eigenvalue.

The following three lemmas are needed for the proofs. Lemma S2.1 and Lemma S2.2 are from [4]; Lemma S2.3 is from [2].

Lemma S2.1. *Suppose $\mathbf{D}$ is an $n \times n$ symmetric matrix whose eigenvalues are all positive except for one negative. If $\mathbf{1}_n\mathbf{1}_n^T + \mathbf{D} \geq 0$, there are $a > 0$, $b > 0$ and $0 \leq s \leq 1$ such that*

$$s\mathbf{1}_n\mathbf{1}_n^T + a\mathbf{I}_n \leq \mathbf{1}_n\mathbf{1}_n^T + \mathbf{D} \leq \mathbf{1}_n\mathbf{1}_n^T + b\mathbf{I}_n.$$

Lemma S2.2. *For an $n \times n$ matrix $\mathbf{A} > \mathbf{0}$ and an $n \times p$ full column rank matrix $\mathbf{h}(\mathbf{x}^{\mathscr{D}})$,*

$$\mathbf{1}_n^T\{\mathbf{A}^{-1} - \mathbf{A}^{-1}\mathbf{h}(\mathbf{x}^{\mathscr{D}})\left(\mathbf{h}^T(\mathbf{x}^{\mathscr{D}})\mathbf{A}^{-1}\mathbf{h}(\mathbf{x}^{\mathscr{D}})\right)^{-1}\mathbf{h}^T(\mathbf{x}^{\mathscr{D}})\mathbf{A}^{-1}\}\mathbf{1}_n = 0$$

*if and only if $\mathbf{1}_n \in \mathcal{C}(\mathbf{h}(\mathbf{x}^{\mathscr{D}}))$.*

Define $\mathbf{A} = \mathbf{R} - \mathbf{1}_n\mathbf{1}_n^T$ for the following lemma and the rest of the proofs in the supplementary materials.

LEMMA S2.3. *If* $\mathbf{1}_n \in \mathcal{C}(\mathbf{h}(\mathbf{x}^{\mathscr{D}}))$*, then* $\mathbf{A}$ *is nonsingular and*

$$\mathbf{R}^{-1}\mathbf{P}_{\mathbf{R}} = \mathbf{A}^{-1}\mathbf{P}_{\mathbf{A}},$$

*where* $\mathbf{P}_{\mathbf{R}}$ *is defined in Equation (2.6) and*

$$\mathbf{P}_{\mathbf{A}} = \mathbf{I}_{\boldsymbol{\gamma}} - \mathbf{h}(\mathbf{x}^{\mathscr{D}})\{\mathbf{h}^T(\mathbf{x}^{\mathscr{D}})\mathbf{A}^{-1}\mathbf{h}(\mathbf{x}^{\mathscr{D}})\}^{-1}\mathbf{h}^T(\mathbf{x}^{\mathscr{D}})\mathbf{A}^{-1}.$$

PROOF OF LEMMA 3.3. (i) If $\forall l$, $1 \le l \le p$, $\gamma_l \to 0^+$, one has $\mathbf{R} \to \mathbf{I}_n$, and the marginal likelihood becomes

$$L(\boldsymbol{\gamma} \mid \mathbf{y}^{\mathscr{D}}) \propto |\mathbf{h}^T(\mathbf{x}^{\mathscr{D}})\mathbf{h}(\mathbf{x}^{\mathscr{D}})|^{-\frac{1}{2}}(S_0^2)^{-(\frac{n-q}{2}+a-1)},$$

where

$$S_0^2 = (\mathbf{y}^{\mathscr{D}} - \mathbf{h}(\mathbf{x}^{\mathscr{D}})\hat{\boldsymbol{\theta}}_0)^T(\mathbf{y}^{\mathscr{D}} - \mathbf{h}(\mathbf{x}^{\mathscr{D}})\hat{\boldsymbol{\theta}}_0),$$

with $\hat{\boldsymbol{\theta}}_0 = (\mathbf{h}^T(\mathbf{x}^{\mathscr{D}})\mathbf{h}(\mathbf{x}^{\mathscr{D}}))^{-1}\mathbf{h}^T(\mathbf{x}^{\mathscr{D}})\mathbf{y}^{\mathscr{D}}$. Similarly, the profile likelihood will be

$$L(\boldsymbol{\gamma} \mid \mathbf{y}^{\mathscr{D}}, \hat{\sigma}^2_{MLE}, \hat{\boldsymbol{\theta}}_{MLE}) \propto (S_0^2)^{-n/2}.$$

Hence the marginal likelihood and profile likelihood exist, and their values are positive.

(ii) Because

$$\mathbf{R} = \mathbf{R}_1 \circ \mathbf{R}_2 \circ \cdots \circ \mathbf{R}_p,$$

and for each $\mathbf{R}_l$, $l = 1, \ldots, p$, $\mathbf{R}_l = \mathbf{1}_n\mathbf{1}_n^T + \nu_l(\gamma_l)(\mathbf{D}_l + o(1))$, it follows from Lemma S2.1 that

$$C_{l1}\mathbf{1}_n\mathbf{1}_n^T + C_{l2}\nu_l(\gamma_l)\mathbf{I}_n \le \mathbf{R}_l \le \mathbf{1}_n\mathbf{1}_n^T + C_{l3}\nu_l(\gamma_l)\mathbf{I}_n,$$

and, from Fact 2, that

$$b_1\mathbf{1}_n\mathbf{1}_n^T + b_2\mathbf{I}_n \le \mathbf{R} \le \mathbf{1}_n\mathbf{1}_n^T + b_3\mathbf{I}_n, \tag{S12}$$

where $b_1 = \prod_l^p C_{l1}$, $b_2 = \prod_l^p\{C_{l1} + C_{l2}\nu_l(\gamma_l)\} - \prod_l^p C_{l1}$ and $b_3 = \prod_l^p\{1 + C_{l3}\nu_l(\gamma_l)\} - 1$.

Using Equation (S10) yields

$$b_2^{n-1}(b_2 + b_1 n) \le \quad |\mathbf{R}| \quad \le b_3^{n-1}(b_3 + n).$$

Thus

$$|\mathbf{R}| = O(b_2^{n-1}) = O\left((\sum_{l=1}^{p}\nu_l(\gamma_l))^{n-1}\right). \tag{S13}$$

Using Equation (S11), it follows that

$$b_3^{-1}(\mathbf{I}_n - \frac{\mathbf{1}_n\mathbf{1}_n^T}{b_3+n}) \le \mathbf{R}^{-1} \le b_2^{-1}(\mathbf{I}_n - \frac{b_1\mathbf{1}_n\mathbf{1}_n^T}{b_2+nb_1}), \tag{S14}$$

and using Equation (S10), that

$$\begin{aligned} & b_3^{-q}\left|\mathbf{h}^T(\mathbf{x}^{\mathscr{D}})\mathbf{h}(\mathbf{x}^{\mathscr{D}})\right|(1-\frac{\mathbf{1}_n^T\mathbf{P}_x\mathbf{1}_n}{b_3+n}) \\ \le\ & \left|\mathbf{h}^T(\mathbf{x}^{\mathscr{D}})\mathbf{R}^{-1}\mathbf{h}(\mathbf{x}^{\mathscr{D}})\right| \le b_2^{-q}\left|\mathbf{h}^T(\mathbf{x}^{\mathscr{D}})\mathbf{h}(\mathbf{x}^{\mathscr{D}})\right|(1-b_1\frac{\mathbf{1}_n^T\mathbf{P}_x\mathbf{1}_n}{b_2+nb_1}), \end{aligned} \tag{S15}$$

where $\mathbf{P}_x = \mathbf{h}(\mathbf{x}^{\mathscr{D}})\left(\mathbf{h}^T(\mathbf{x}^{\mathscr{D}})\mathbf{h}(\mathbf{x}^{\mathscr{D}})\right)^{-1}\mathbf{h}^T(\mathbf{x}^{\mathscr{D}})$. Thus if $\mathbf{1}_n \notin \mathcal{C}(\mathbf{h}(\mathbf{x}^{\mathscr{D}}))$,

$$\left|\mathbf{h}^T(\mathbf{x}^{\mathscr{D}})\mathbf{R}^{-1}\mathbf{h}(\mathbf{x}^{\mathscr{D}})\right| = O\left((\sum_{l=1}^{p}\nu_l(\gamma_l))^{-q}\right). \tag{S16}$$

If $\mathbf{1}_n \in \mathcal{C}(\mathbf{h}(\mathbf{x}))$, by applying Lemma S2.2, one additionally has $\mathbf{1}_n^T\mathbf{P}_x\mathbf{1}_n = n$. Applying this fact to Inequality (S15) yields

$$\begin{aligned} & b_3^{-(q-1)}\left|\mathbf{h}^T(\mathbf{x}^{\mathscr{D}})\mathbf{h}(\mathbf{x}^{\mathscr{D}})\right|\frac{1}{b_3+n} \\ \le\ & \left|\mathbf{h}^T(\mathbf{x}^{\mathscr{D}})\mathbf{R}^{-1}\mathbf{h}(\mathbf{x}^{\mathscr{D}})\right| \le b_2^{-(q-1)}\left|\mathbf{h}^T(\mathbf{x}^{\mathscr{D}})\mathbf{h}(\mathbf{x}^{\mathscr{D}})\right|\frac{1}{b_2+nb_1}, \end{aligned}$$

from which it follows that

$$\left|\mathbf{h}^T(\mathbf{x}^{\mathscr{D}})\mathbf{R}^{-1}\mathbf{h}(\mathbf{x}^{\mathscr{D}})\right| = O\left((\sum_{l}^{p}\nu_l(\gamma_l))^{-(q-1)}\right). \tag{S17}$$

According to Equation (S14), $\mathbf{R}^{-1} = O(b_4^{-1}(\mathbf{I}_n - \frac{\mathbf{1}_n\mathbf{1}_n^T}{b_5+n}))$, where $b_4^{-1} = O((\sum_l^p(\nu_l(\gamma_l)))) \to 0$ and $b_5^{-1} = O((\sum_l^p(\nu_l(\gamma_l)))) \to 0$, when $\gamma_l \to \infty$ for all $1 \le l \le p$. Plugging $\mathbf{R}^{-1} = O(b_4^{-1}(\mathbf{I}_n - \frac{\mathbf{1}_n\mathbf{1}_n^T}{b_5+n}))$ into $\mathbf{Q}$ below, if $\mathbf{1}_n \notin \mathcal{C}(\mathbf{h}(\mathbf{x}^{\mathscr{D}}))$, then

$$\begin{aligned} \mathbf{Q} &= \mathbf{R}^{-1} - \mathbf{R}^{-1}\mathbf{h}(\mathbf{x}^{\mathscr{D}})\{\mathbf{h}^T(\mathbf{x}^{\mathscr{D}})\mathbf{R}^{-1}\mathbf{h}(\mathbf{x}^{\mathscr{D}})\}^{-1}\mathbf{h}^T(\mathbf{x}^{\mathscr{D}})\mathbf{R}^{-1} \\ &= O\left(b_4^{-1}(\mathbf{I}_n - \mathbf{P}_x - \frac{(\mathbf{I}_n-\mathbf{P}_x)\mathbf{1}_n\mathbf{1}_n^T(\mathbf{I}_n-\mathbf{P}_x)}{b_5+n-\mathbf{1}_n^T\mathbf{P}_x\mathbf{1}_n})\right) \\ &= O\left((\sum_{l=1}^{p}\nu_l(\gamma_l))^{-1}\left(\mathbf{I}_n - \mathbf{P}_x - \frac{\mathbf{1}_n\mathbf{1}_n^T}{n}\right)\right). \end{aligned} \tag{S18}$$

If $\mathbf{1}_n \in \mathcal{C}(\mathbf{h}(\mathbf{x}^{\mathscr{D}}))$, using Equation (S12) and the fact that $\mathbf{1}_n\mathbf{1}_n^T$ is positive semidefinite, one has

$$b_3^{-1}\mathbf{I}_n \leq \mathbf{A}^{-1} \leq b_2^{-1}\left(\mathbf{I}_n - \frac{b_1 - 1}{b_2 + n(b_1 - 1)}\mathbf{1}_n\mathbf{1}_n^T\right) \leq b_2^{-1}\mathbf{I}_n,$$

where $\mathbf{A}$ is defined before Lemma S2.3. Define $b_6^{-1} = O((\sum_l^p(\nu_l(\gamma_l)))) \to 0$, when $\gamma_l \to \infty$ for all $1 \leq l \leq p$. Using Lemma S2.3, it follows that

$$\begin{aligned}\mathbf{Q} &= \mathbf{A}^{-1} - \mathbf{A}^{-1}\mathbf{h}(\mathbf{x}^{\mathscr{D}})\{\mathbf{h}^T(\mathbf{x}^{\mathscr{D}})\mathbf{A}^{-1}\mathbf{h}(\mathbf{x}^{\mathscr{D}})\}^{-1}\mathbf{h}^T(\mathbf{x}^{\mathscr{D}})\mathbf{A}^{-1} \\ &= O\left(b_6^{-1}(\mathbf{I} - \mathbf{P}_x)\right) \\ &= O\left((\sum_l^p \nu_l(\gamma_l))^{-1}(\mathbf{I} - \mathbf{P}_x)\right).\end{aligned} \tag{S19}$$

Using Equation (S18) and Equation (S19), we have

$$S^2 = (y^{\mathscr{D}})^T\mathbf{Q}\mathbf{y}^{\mathscr{D}} = O\left((\sum_l^p \nu_l(\gamma_l))^{-1}\right). \tag{S20}$$

By combining Equation (S13), (S16), (S17) and (S20), the proof is complete.

□

Proof of Lemma 3.4. (i) If $\forall l$, $1 \leq l \leq p$, $\gamma_l \to 0^+$, one has $\mathbf{R} \to \mathbf{I}_n$ and

$$\mathbf{Q} \to \mathbf{I}_n - \mathbf{P}_x.$$

For $\forall l$, $1 \leq l \leq p$, we have

$$\mathrm{tr}(\mathbf{W}_l^2) = \mathrm{tr}\left[\left(\frac{\partial \mathbf{R}}{\partial \gamma_l}\mathbf{Q}\right)^2\right] \leq C\mathrm{tr}\left[\left(\frac{\partial \mathbf{R}}{\partial \gamma_l}\right)^2\right], \tag{S21}$$

with $C > 0$ being a constant. Note that $\mathbf{I}^*(\boldsymbol{\gamma})$ is the Fisher information matrix from the marginal likelihood and is positive semidefinite. According to the Hadamard's inequality, the determinant of a positive semidefinite matrix is bounded by its diagonal elements, and thus,

$$\pi^R(\boldsymbol{\gamma}) \propto |\mathbf{I}^*(\boldsymbol{\gamma})|^{1/2} \leq \left[(n-q)\prod_{l=1}^p \mathrm{tr}(\mathbf{W}_l^2)\right]^{1/2} \leq C\left[\prod_{l=1}^p \mathrm{tr}\left(\frac{\partial \mathbf{R}}{\partial \gamma_l}\right)^2\right]^{1/2}. \tag{S22}$$

with $C > 0$.

(ii) For simplicity, $\nu_l$ and $\omega_l$ are used to represent $\nu_l(\gamma_l)$ and $\omega_l(\gamma_l)$ in the proof, respectively. If $\mathbf{1}_n \notin \mathcal{C}(\mathbf{h}(\mathbf{x}^{\mathscr{D}}))$, Assumption 3.2 implies, for any $1 \le m \le p$,

$$\mathbf{R} = \mathbf{1}_n\mathbf{1}_n^T + \sum_{l=1}^{p} \nu_l \mathbf{D}_l + \mathbf{F}_{-m} + o(\nu_m),$$

where $\mathbf{F}_{-m}$ is an $n \times n$ matrix that does not depend on $\gamma_m$. Thus, for any $1 \le l \le p$,

$$\left\| \frac{\partial \mathbf{R}}{\partial \gamma_l} \right\|_\infty \le C|\nu_l'|, \tag{S23}$$

where $C > 0$ is a positive constant. Using Equation (S18) yields

$$\mathrm{tr}(\mathbf{W}_l^2) = \mathrm{tr}\left(\frac{\partial \mathbf{R}}{\partial \gamma_l}\mathbf{Q}\right)^2 = O\left(\left(\frac{\nu_l'}{\sum_{l=1}^p \nu_l}\right)^2\right), \tag{S24}$$

and thus Equation (3.2) follows using the fact that the determinant of positive semidefinite matrix is bounded by its diagonal elements.
For the case $\mathbf{1}_n \in \mathcal{C}(\mathbf{h}(\mathbf{x}^{\mathscr{D}}))$, w.l.o.g., we assume $m = 1$. Denote

$$\mathbf{\Psi} = \sum_{l=1}^{p} \frac{\nu_l}{\nu_l'} \frac{\partial \mathbf{R}}{\partial \gamma_l} - \mathbf{A}, \tag{S25}$$

and

$$\frac{\partial \mathbf{R}_1}{\partial \gamma_1} = \frac{\nu_1'}{\nu_1}\left(\mathbf{\Psi} + \mathbf{A} - \sum_{l=2}^{p} \frac{\nu_l}{\nu_l'} \frac{\partial \mathbf{R}}{\partial \gamma_l}\right). \tag{S26}$$

Using Lemma S2.3 yields

$$\begin{aligned} \mathrm{tr}(\mathbf{W}_1) &= \frac{\nu_1'}{\nu_1}\mathrm{tr}\left[\left(\mathbf{\Psi} + \mathbf{A} - \sum_{l=2}^{p} \frac{\nu_l}{\nu_l'} \frac{\partial \mathbf{R}}{\partial \gamma_l}\right)\mathbf{A}^{-1}\mathbf{P}_{\mathbf{A}}\right] \\ &= \frac{\nu_1'}{\nu_1}\mathrm{tr}\left[\mathbf{\Psi}\mathbf{A}^{-1}\mathbf{P}_{\mathbf{A}} + \mathbf{P}_{\mathbf{A}} - \sum_{l=2}^{p} \frac{\nu_l}{\nu_l'}\mathbf{W}_l\right], \end{aligned} \tag{S27}$$

$$\mathrm{tr}(\mathbf{W}_1^2) = \left(\frac{\nu_1'}{\nu_1}\right)^2 \mathrm{tr}\left[\mathbf{\Psi}\mathbf{A}^{-1}\mathbf{P}_{\mathbf{A}} + \mathbf{P}_{\mathbf{A}} - \sum_{l=2}^{p} \frac{\nu_l}{\nu_l'}\mathbf{W}_l\right]^2, \tag{S28}$$

and, for $2 \leq j \leq p$,

$$\text{tr}(\mathbf{W}_1\mathbf{W}_j) = \frac{\nu_1'}{\nu_1}\text{tr}\left[\left(\mathbf{\Psi}\mathbf{A}^{-1}\mathbf{P}_{\mathbf{A}} + \mathbf{P}_{\mathbf{A}} - \sum_{l=2}^{p}\frac{\nu_l}{\nu_l'}\mathbf{W}_l\right)\mathbf{W}_j\right]. \tag{S29}$$

Note that $(\mathbf{P}_{\mathbf{A}})^2 = \mathbf{P}_{\mathbf{A}}$, $\text{tr}(\mathbf{W}_l\mathbf{P}_{\mathbf{A}}) = \text{tr}(\mathbf{W}_l)$ for all $1 \leq l \leq p$, and $\text{tr}(\mathbf{P}_{\mathbf{A}}) = n - q$. For the reference prior defined in Equation (2.8), first put $\frac{\nu_1'}{\nu_1}$ outside the determinant by dividing $\frac{\nu_1'}{\nu_1}$ on the second column and the second row. Then multiply the first row by $-1$ and add it to the second row. Also, multiply the first column by $-1$ and add it to the second column. After the above manipulation of the determinant, it follows that

(S30)

$$\pi^R(\boldsymbol{\gamma}) \propto \frac{\nu_1'}{\nu_1}\begin{vmatrix} n-q & \text{tr}\left(\mathbf{B}\right) & \text{tr}(\mathbf{W}_2) & \dots & \text{tr}(\mathbf{W}_p) \\ & \text{tr}\left(\mathbf{B}^2\right) & \text{tr}\left(\mathbf{B}\mathbf{W}_2\right) & \dots & \text{tr}\left(\mathbf{B}\mathbf{W}_p\right) \\ & & \text{tr}(\mathbf{W}_2^2) & \dots & \text{tr}(\mathbf{W}_2\mathbf{W}_p) \\ & & & \ddots & \vdots \\ & & & & \text{tr}(\mathbf{W}_p^2) \end{vmatrix}^{1/2},$$

where $\mathbf{B} = \mathbf{\Psi}\mathbf{A}^{-1}\mathbf{P}_{\mathbf{A}} - \sum_{l=2}^{p}\frac{\nu_l}{\nu_l'}\mathbf{W}_l$. Further multiple the $(l+1)^{th}$ column by $\frac{\nu_1}{\nu_1'}$, $2 \leq l \leq p$, and add them to the second column. Then multiply the $(l+1)^{th}$ row by $\frac{\nu_1}{\nu_1'}$, $2 \leq l \leq p$, and add them to the second row. It follows that

(S31)

$$\pi^R(\boldsymbol{\gamma}) \propto \frac{\nu_1'}{\nu_1}\begin{vmatrix} n-p & \text{tr}\left(\mathbf{\Psi}\mathbf{A}^{-1}\mathbf{P}_{\mathbf{A}}\right) & \text{tr}(\mathbf{W}_2) & \dots & \text{tr}(\mathbf{W}_p) \\ & \text{tr}\left(\mathbf{\Psi}\mathbf{A}^{-1}\mathbf{P}_{\mathbf{A}}\right)^2 & \text{tr}\left\{\mathbf{\Psi}\mathbf{A}^{-1}\mathbf{P}_{\mathbf{A}}\mathbf{W}_2\right\} & \dots & \text{tr}\left\{\mathbf{\Psi}\mathbf{A}^{-1}\mathbf{P}_{\mathbf{A}}\mathbf{W}_p\right\} \\ & & \text{tr}(\mathbf{W}_2^2) & \dots & \text{tr}(\mathbf{W}_2\mathbf{W}_p) \\ & & & \ddots & \vdots \\ & & & & \text{tr}(\mathbf{W}_p^2) \end{vmatrix}^{1/2}.$$

By definition,

$$\mathbf{A} = \sum_{l=1}^{p}\nu_l\mathbf{D}_l + \sum_{l=1}^{p}\nu_l\omega_l\mathbf{D}_l^* + \sum_{l\neq m}\nu_l\nu_m(\mathbf{D}_l \circ \mathbf{D}_m + o(1)), \tag{S32}$$

and

(S33)

$$\frac{\partial \mathbf{R}}{\partial \gamma_l} = \nu_l'\mathbf{D}_l + [\nu_l'\omega_l + \nu_l\omega_l']\mathbf{D}_l^* + \nu_l'\sum_{m\in\{1,\dots,p\}\backslash l}\nu_m(\mathbf{D}_l \circ \mathbf{D}_m + o(1)).$$

By Equation (S25), we have

$$\boldsymbol{\Psi} = O(\sum_{l=1}^{p} \frac{\nu_l^2 \omega_l'}{\nu_l'} \mathbf{D}_l^*). \tag{S34}$$

Directly applying Cramer's rule to the reference prior in Equation (S31) yields

$$O\left(\frac{\nu_1'}{\nu_1}\{\prod_{l=2}^{p} \mathrm{tr}(\mathbf{W}_l^2) \mathrm{tr}\left(\boldsymbol{\Psi}\mathbf{A}^{-1}\mathbf{P}_{\mathbf{A}}\right)^2\}^{1/2}\right).$$

Using Equation (S19), (S24) and (S34), the result for $p \geq 2$ of (ii) in Lemma 3.4 follows.

□

Proof of Lemma 3.5. For any $\gamma_{l_1} \to \infty$, $1 \leq l_1 \leq p_1$, $\mathbf{R}_{l_1} \to \mathbf{1}_n \mathbf{1}_n^T$ and the correlation matrix is $\mathbf{R} \to \mathbf{R}_{p_1+1} \circ \mathbf{R}_{p_1+2} \circ ... \circ \mathbf{R}_p$. Note that $\gamma_l < \infty$ with $p_1 + 1 \leq l \leq p$. The following result is from the first part of Lemma 3.3,

$$L(\boldsymbol{\gamma} \mid \mathbf{y}^{\mathscr{D}}) = O(1) \tag{S35}$$

Since at least one $\gamma_l$ does not go to infinity, Equation (S18) and Equation (S19) yield,

$$\mathbf{Q} = O(1).$$

Assumption 3.2 implies

$$\left\| \frac{\partial \mathbf{R}}{\partial \gamma_l} \right\|_{\infty} \leq C |\nu_l'(\gamma_l)|, \tag{S36}$$

where $C$ is a postive constant. So that, for any $1 \leq l_1 \leq p_1$,

$$\mathrm{tr}(\mathbf{W}_{l_1}^2) = \mathrm{tr}\left(\left(\frac{\partial \mathbf{R}}{\partial \gamma_{l_1}}\mathbf{Q}\right)^2\right) = O\left(\nu_{l_1}'(\gamma_{l_1})^2\right). \tag{S37}$$

Since $p_1 + 1 \leq l_2 \leq p_2$, $\gamma_l \to 0^+$, the proof of Lemma 3.4 yields

$$\mathrm{tr}(\mathbf{W}_{l_2}^2) \leq C \mathrm{tr}\left(\frac{\partial \mathbf{R}}{\partial \gamma_{l_2}}\right)^2.$$

And for $p_1 + 1 \leq l_3 \leq p$, since $\gamma_{l_3}$ is finite,

$$\mathrm{tr}(\mathbf{W}_{l_3}^2) = \mathrm{tr}\left(\frac{\partial \mathbf{R}}{\partial \gamma_{l_3}}\mathbf{Q}\right)^2 = O(1). \tag{S38}$$

The fact that the determinant of a positive semidefinite matrix is bounded by its diagonal elements yields

$$\pi^R(\boldsymbol{\gamma}) \leq C \left| \prod_{l_1=1}^{p_1} \nu'_{l_1}(\gamma_{l_1}) \left[ \prod_{l_2=p_1+1}^{p_2} \mathrm{tr}\left( \frac{\partial \mathbf{R}}{\partial \gamma_{l_2}} \right)^2 \right]^{1/2} \right|. \tag{S39}$$

Combining the result of Equation (S35) and Equation (S39), the proof is complete.

□

PROOF OF THEOREM 3.2. Only the proof for Matérn correlation is given here; and the rest of cases can be checked similarly. When $p = 1$, posterior propriety is established in [2]. For $p > 1$, only the case $\mathbf{1}_n \in \mathcal{C}(\mathbf{h}(\mathbf{x}^{\mathscr{D}}))$ is shown below; the case $\mathbf{1}_n \notin \mathcal{C}(\mathbf{h}(\mathbf{x}^{\mathscr{D}}))$ can be checked similarly.

(i) First assume $\gamma_{(1)} \leq \gamma_{(2)} \leq ... \leq \gamma_{(p)}$ and each $\gamma_l$ goes to $\infty$ .

(i.1) For the case $0 < \alpha < 1$, the marginal posterior is

$$\pi^R(\boldsymbol{\gamma}) L(\boldsymbol{\gamma} \mid \mathbf{y}^{\mathscr{D}}) \leq C_1 \frac{\prod_{l=1}^{p} \gamma_l^{-2\alpha-1}}{\gamma_{(1)}^{-2p\alpha}} \gamma_{(1)}^{-2+2\alpha}, \tag{S40}$$

with $C_1 > 0$ being a constant. To show that Equation (S40) is integrable, we need only prove that $\int_M^\infty \int_{\gamma_{(1)}}^\infty ... \int_{\gamma_{(p-1)}}^\infty \pi^R(\boldsymbol{\gamma}) L(\boldsymbol{\gamma} \mid \mathbf{y}^{\mathscr{D}}) d\boldsymbol{\gamma}$ is finite, which is easily seen from the following:

$$\begin{aligned}
&\int_M^\infty \int_{\gamma_{(1)}}^\infty ... \int_{\gamma_{(p-1)}}^\infty \pi^R(\boldsymbol{\gamma}) L(\boldsymbol{\gamma} \mid \mathbf{y}^{\mathscr{D}}) d\boldsymbol{\gamma} \\
\leq\ & C_2 \int_M^\infty \int_{\gamma_{(1)}}^\infty ... \int_{\gamma_{(p-1)}}^\infty \prod_{l=2}^{p} \gamma_{(l)}^{-2\alpha-1} \gamma_{(1)}^{2p\alpha-3} d\gamma_{(p)} ... d\gamma_{(1)} \\
=\ & C_3 \int_M^\infty \gamma_{(1)}^{2\alpha-3} d\gamma_{(1)} \\
=\ & C_4 M^{2\alpha-2} \\
<\ & \infty,
\end{aligned}$$

for $M > 0$, $C_2 > 0$, $C_3 > 0$, $C_4 > 0$ being constants and $0 < \alpha < 1$. By Fubini's Theorem, Equation (S40) is integrable.

(i.2) For the case $\alpha = 1$, as

$$\pi^R(\boldsymbol{\gamma})L(\boldsymbol{\gamma} \mid \mathbf{y}^{\mathscr{D}}) \leq C_1 \frac{\prod_{l=1}^{p} \frac{2\log\gamma_l - 1}{\gamma_l^3}}{\left(\frac{\log(\gamma_1)}{\gamma_{(1)}^2}\right)^p} \frac{\frac{1}{\gamma_{(1)}^2(2\log\gamma_{(1)}-1)}}{\frac{\log(\gamma_{(1)})}{\gamma_{(1)}^2}},$$

we have

$$\begin{aligned}
&\int_M^\infty \int_{\gamma_{(1)}}^\infty \cdots \int_{\gamma_{(p-1)}}^\infty \pi^R(\boldsymbol{\gamma})L(\boldsymbol{\gamma} \mid \mathbf{y}^{\mathscr{D}})d\boldsymbol{\gamma} \\
\leq\ & \int_M^\infty \int_{\gamma_{(1)}}^\infty \cdots \int_{\gamma_{(1)}}^\infty \pi^R(\boldsymbol{\gamma})L(\boldsymbol{\gamma} \mid \mathbf{y}^{\mathscr{D}})d\boldsymbol{\gamma} \\
\leq\ & C_2 \int_M^\infty \frac{\frac{2\log\gamma_{(1)}-1}{\gamma_{(1)}^3}\left(\frac{\log(\gamma_{(1)})}{\gamma_{(1)}^2}\right)^{p-1}}{\left(\frac{\log\gamma_{(1)}}{\gamma_{(1)}^2}\right)^{p+1}\gamma_{(1)}^2(2\log\gamma_{(1)}-1)} d\gamma_{(1)} \\
=\ & C_3 \int_M^\infty \frac{1}{\gamma_{(1)}\log^2\gamma_{(1)}} d\gamma_{(1)} \\
=\ & \frac{C_4}{\log M} \\
<\ & \infty.
\end{aligned}$$

for $C_1 > 0$, $C_2 > 0$, $C_3 > 0$, $C_4 > 0$ and $M > 0$ being constants.

(i.3) For the case $1 < \alpha < 2$, similarly we have,

$$\begin{aligned}
&\int_M^\infty \int_{\gamma_{(1)}}^\infty \cdots \int_{\gamma_{(p-1)}}^\infty \pi^R(\boldsymbol{\gamma})L(\boldsymbol{\gamma} \mid \mathbf{y}^{\mathscr{D}})d\boldsymbol{\gamma} \\
\leq\ & C_1 \int_M^\infty \int_{\gamma_{(1)}}^\infty \cdots \int_{\gamma_{(p-1)}}^\infty \frac{\prod_{l=1}^{p}\gamma_l^{-3}}{\gamma_{(1)}^{-2p}} \gamma_{(1)}^{2-2\alpha} d\gamma_{(p)} \ldots d\gamma_{(1)} \\
=\ & C_2 \int_M^\infty \gamma_{(1)}^{-2\alpha+1} d\gamma_{(1)} \\
=\ & C_3 M^{-2\alpha+2} \\
<\ & \infty.
\end{aligned}$$

for $C_1 > 0$, $C_2 > 0$, $C_3 > 0$ and $M > 0$ being constants.

(i.4) For the case $\alpha = 2$, we have

$$
\begin{aligned}
&\int_M^\infty \int_{\gamma_{(1)}}^\infty \dots \int_{\gamma_{(p-1)}}^\infty \pi^R(\boldsymbol{\gamma}) L(\boldsymbol{\gamma} \mid \mathbf{y}^{\mathscr{D}}) d\boldsymbol{\gamma} \\
\leq\ & C_1 \int_M^\infty \int_{\gamma_{(1)}}^\infty \dots \int_{\gamma_{(p-1)}}^\infty \frac{\prod_{l=1}^p \gamma_l^{-3}}{\gamma_{(1)}^{-2p}} \frac{(2\log(\gamma_{(1)}) - 1)}{\gamma_{(1)}^2} d\gamma_{(p)} \dots d\gamma_{(1)} \\
=\ & C_2 \int_M^\infty \gamma_{(1)}^{-3} (2\log(\gamma_{(1)}) - 1) d\gamma_{(1)} \\
=\ & C_3 \frac{\log M}{M^2} \\
<\ & \infty.
\end{aligned}
$$

for $C_1 > 0$, $C_2 > 0$, $C_3 > 0$ and $M > 0$ being constants.

(i.5) For the case $\alpha > 2$,

$$
\pi^R(\boldsymbol{\gamma}) L(\boldsymbol{y}^{\mathscr{D}} | \boldsymbol{\gamma}) \leq C \frac{\prod_{l=1}^p \gamma_l^{-3}}{\gamma_{(1)}^{-2p}} \gamma_{(1)}^{-2} \leq C \prod_{l=1}^p \gamma_l^{-1-2/p}.
$$

for $C > 0$ being the constant. The right hand side is clearly integrable.

(ii) If there is at least one $l$ for which $\gamma_l < \infty$, Lemma 3.5 shows that the integral of the product is just the product of the individual integrals, so one only needs to check that $\nu'(\gamma_l)$ is integrable when $\gamma_l \to \infty$ and that $\frac{\partial R}{\partial \gamma_l}$ is integrable when $\gamma_l \to 0$. From Table 1, $\nu'(\gamma_l)$ is integrable when $\gamma_l \to \infty$. Noting that by the property of the modified Bessel function of the second kind, when $z \to \infty$ ([1], Section 9.7.4.)

$$
\frac{\partial \mathcal{K}_\alpha(z)}{\partial z} \to -\sqrt{\frac{\pi}{2z}} \exp(-z)(1 + O(\frac{1}{z})),
$$

from which the Matérn correlation function is integrable as $\gamma_l \to 0$.

□

## S3. Proof for Section 4.3.

Proof of Lemma 4.1. Since

$$
\tilde{\mathbf{R}} = \mathbf{R} + \eta \mathbf{I}_n
$$

and applying similar derivations as in the proof of Lemma 3.3, it is easy to see that

$$|\tilde{\mathbf{R}}| = O\left( (\sum_{l=1}^{p} \nu_l(\gamma_l) + \eta)^{n-1} \right) \; : \tag{S41}$$

(S42)

$$\left| \mathbf{h}^T(\mathbf{x}^{\mathscr{D}})\tilde{\mathbf{R}}^{-1}\mathbf{h}(\mathbf{x}^{\mathscr{D}}) \right| = \begin{cases} O\left( (\sum_{l=1}^{p} \nu_l(\gamma_l) + \eta)^{-q} \right), & \mathbf{1}_n \notin \mathcal{C}(\mathbf{h}(\mathbf{x}^{\mathscr{D}})), \\ O\left( (\sum_{l=1}^{p} \nu_l(\gamma_l) + \eta)^{-(q-1)} \right), & \mathbf{1}_n \in \mathcal{C}(\mathbf{h}(\mathbf{x}^{\mathscr{D}})) \; : \end{cases}$$

$$\tilde{\mathbf{S}}^2 = (\mathbf{y}^{\mathscr{D}})^T \tilde{\mathbf{Q}} \mathbf{y}^{\mathscr{D}} = O\left( (\sum_{l=1}^{p} \nu_l(\gamma_i) + \eta)^{-1} \right). \tag{S43}$$

The result then immediately follows. □

Proof of Lemma 4.2. The proofs of (a) and (b) can be done similarly to the Proof of Lemma 3.4, by noting that $\tilde{\mathbf{R}} = \mathbf{R}_1 \circ \mathbf{R}_2 \circ ... \circ \mathbf{R}_p \circ \mathbf{R}_{p+1}$, where $\mathbf{R}_{p+1} = \eta \mathbf{I}_n + \mathbf{1}_n \mathbf{1}_n^T$, with $\nu_{p+1}(\eta) = \eta$ and $\omega_{p+1}(\eta) = 0$. □

**S4. Figure for Section 5.3.** Figure S1 presents plots of the borehole function made by fixing seven of the inputs and varying one. From Figure S1, the outputs barely change when the $2^{nd}$, $3^{rd}$ and $5^{th}$ inputs (i.e., $r$, $T_u$ and $T_l$ in the borehole function) vary.

## References.

[1] Abramowitz, M., Stegun, I. A. et al. (1966). *Handbook of mathematical functions*. National Bureau of Standards-Applied Mathematics Series.

[2] Berger, J. O., De Oliveira, V. and Sansó, B. (2001). Objective Bayesian analysis of spatially correlated data. *Journal of the American Statistical Association* **96** 1361–1374.

[3] Lopes, D. (2011). Development and implementation of Bayesian computer model emulators PhD thesis, Duke University.

[4] Ren, C., Sun, D. and He, C. (2012). Objective Bayesian analysis for a spatial model with nugget effects. *Journal of Statistical Planning and Inference* **142** 1933–1946.

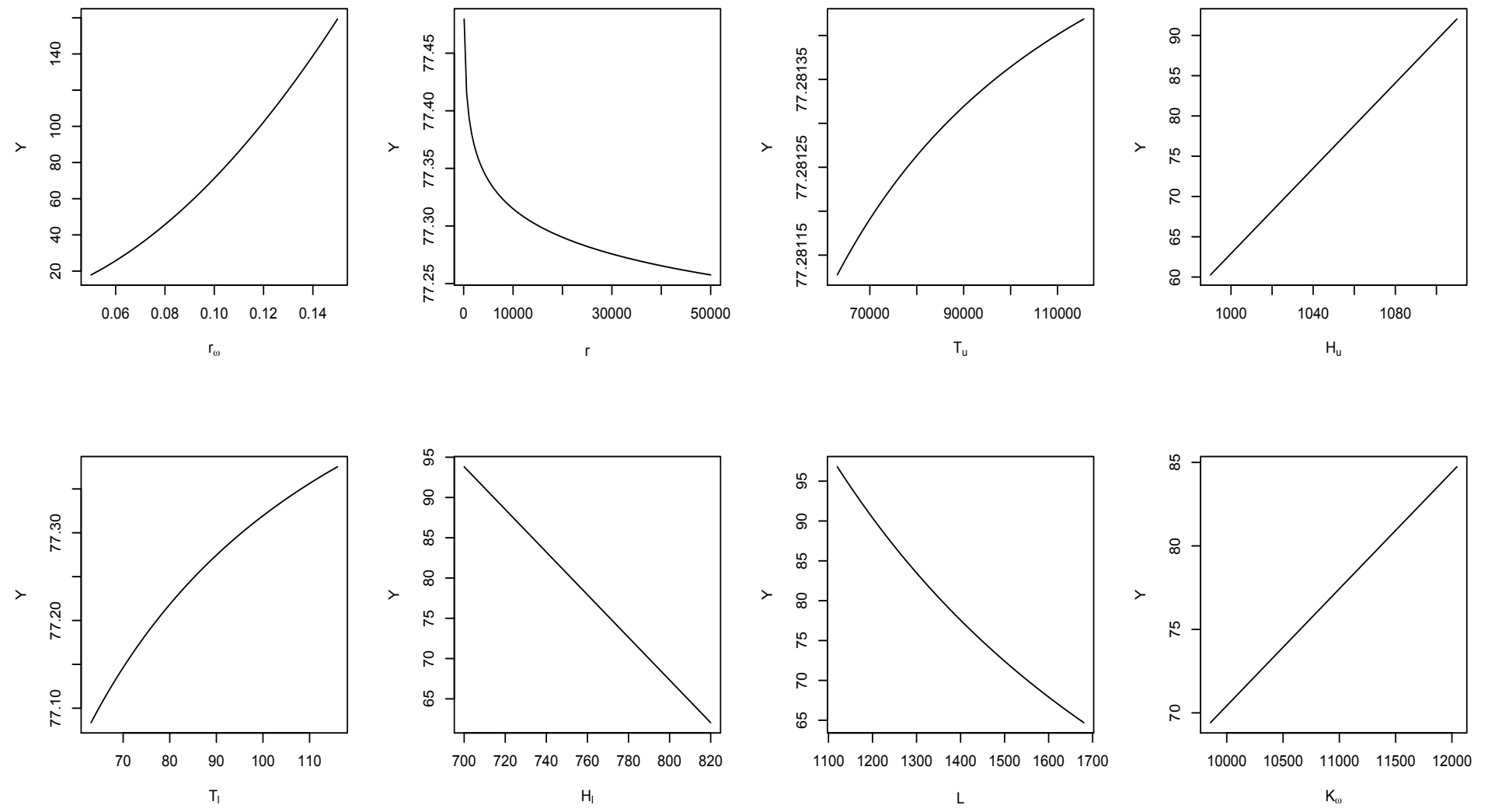

FIG S1. *The plot of borehole output by varying one input at a time.*